\documentclass[reqno,12pt]{amsart}
\usepackage[utf8]{inputenc}
\usepackage{mathtools}
\usepackage{amsmath}
\usepackage{amssymb}
\usepackage{bbm}
\usepackage{tikz-network}
\usepackage{float}
\usepackage{enumitem}
\usepackage{fbb}
\usepackage{microtype}

\usepackage[colorlinks=true,linkcolor=blue,citecolor=blue,pdfpagelabels=false]{hyperref}
\usepackage[nameinlink]{cleveref}

\usepackage{amsthm}
\usepackage{aliascnt}
\usepackage{cleveref}

\theoremstyle{plain}

\newtheorem{theorem}{Theorem}[section]

\newaliascnt{lemma}{theorem}
\newtheorem{lemma}[lemma]{Lemma}
\aliascntresetthe{lemma}

\newaliascnt{proposition}{theorem}

\aliascntresetthe{proposition}

\newaliascnt{corollary}{theorem}
\newtheorem{corollary}[corollary]{Corollary}
\aliascntresetthe{corollary}

\theoremstyle{definition}
\newtheorem{definition}{Definition}[section]

\theoremstyle{remark}
\newtheorem{remark}{Remark}[section]

\crefname{theorem}{theorem}{theorems}
\Crefname{theorem}{Theorem}{Theorems}

\crefname{lemma}{lemma}{lemmas}
\Crefname{lemma}{Lemma}{Lemmas}

\crefname{proposition}{proposition}{propositions}
\Crefname{proposition}{Proposition}{Propositions}

\crefname{corollary}{corollary}{corollaries}
\Crefname{corollary}{Corollary}{Corollaries}

\crefname{definition}{definition}{definitions}
\Crefname{definition}{Definition}{Definitions}

\crefname{remark}{remark}{remarks}
\Crefname{remark}{Remark}{Remarks}

\crefname{example}{example}{examples}
\Crefname{example}{Example}{Examples}

\newcommand{\G}{\overrightarrow{\mathbb{G}}}
\newcommand{\E}{\overrightarrow{E}}

\newcommand{\C}[2]{C_{#1}^{(#2)}}

\newcommand{\us}[2]{u_{#1}^{(#2)}}
\newcommand{\Z}[2]{Z_{#1}^{(#2)}}

\newcommand{\bb}[1]{\mathbbm{#1}}

\newcommand{\Q}{\mathbb{P}}

\newcommand{\Qr}[1]{\mathbb{P}\left(#1\right)}

\newcommand{\e}{\overrightarrow{e}}
\newcommand{\T}[2]{\mathcal{T}_{#1}^{#2}}
\newcommand{\Fo}[1]{\mathbb{F}_{#1}}
\newcommand{\D}[2]{D_{(#1,#2)}}
\newcommand{\W}[2]{W_{(#1,#2)}}
\newcommand{\R}[1]{\mathcal{R}(#1)}
\newcommand{\Ex}[1]{\mathbb{E}\left[#1\right]}
\newcommand{\be}{\begin{equation}}
\newcommand{\ee}{\end{equation}}
\newcommand{\longminus}{\mathbin{\scalebox{2}[1]{$-$}}}
\newcommand{\st}{\stackrel{st}{\le}}
\newcommand{\eqd}{\stackrel{d}{=}}

\title{Interacting Urn Schemes on Infinite Graphs: A Representation Approach}

	\author[A. Bandyopadhyay]{Antar Bandyopadhyay}
    \author[D. Das]{Deborshi Das}
	\address{Theoretical Statistics and Mathematics Unit, Indian Statistical Institute, Delhi, India}
	\email{antar@isid.ac.in and deborshid9@gmail.com}
	
\makeatletter
\@removefromreset{figure}{chapter}
\@addtoreset{figure}{section}
\makeatother

\begin{document}

\begin{abstract}
Interacting urn schemes have emerged as an important framework for studying
reinforcement phenomena in interacting stochastic systems (see, e.g.,
\cite{launay2011interacting,launay2012generalized,WARM2014,aletti2017interacting,Synchronisation-2,kaur2023interacting,kaur2023interacting2}).
Most existing results, however, are restricted to finite systems. In this
article, we investigate interacting urn schemes involving infinitely many urns,
where the interactions are naturally described by a directed graph. The
reinforcement or drawing mechanism at each urn depends on the compositions
of its in-neighboring urns. To analyze the asymptotic behavior of these
systems, we develop a representation of the urn dynamics in terms of
branching Markov processes evolving on a class of random forests introduced
in this article. This extends the representation approach of
\cite{bandyopadhyay2022new}, originally developed for single P\'{o}lya urn
schemes, to interacting urn schemes on locally finite directed graphs.
Using this representation, we establish almost sure convergence of the
color proportions under standard conditions on the reinforcement matrix
or kernel. In particular, we show that the convergence holds independently
of the topology of the underlying graph, and hence our results also apply
to undirected graphs, where the interactions are bidirectional.
We further investigate measure-valued generalizations and other related
reinforced stochastic processes, illustrating the broad applicability of
the representation approach. Finally, we discuss connections with
self-organized criticality.
\end{abstract}

\keywords{Interacting urn schemes; measure-valued urns; reinforced stochastic processes; infinite directed graphs; branching Markov processes, elephant random walks.}
\subjclass[2020]{Primary 60K35; Secondary 60F15, 05C20.}

\maketitle
\pagestyle{plain}

\tableofcontents

\section{Introduction}
Urn models are fundamental probabilistic models that have been studied since the seminal work of Eggenberger and Pólya \cite{eggenberger1923statistik}, where the classical Pólya urn was introduced as a simple example of a reinforced stochastic process. This was followed by Pólya's influential work \cite{polya1930quelques}, which provided an interpretation of the model as a mechanism for contagion and epidemic spread. More than a century later, urn models and their numerous generalizations continue to attract considerable attention. In particular, a wide variety of finite-color urn models and their asymptotic properties have been studied in recent decades (see, e.g., \cite{gouet1989martingale,gouet1993martingale,gouet1997strong,Janson-FCLT,janson2006limit,bose2009strong,dasgupta2011strong,janson2024almost}). Urn models with infinitely many colors have also received considerable attention (see, e.g., \cite{bandyopadhyay2017polya,Measure_valued_Polya_urn,bandyopadhyay2020strong,Infinite_color_urn_Janson,bandyopadhyay2022new,Exchangable_MVPU}). Beyond their intrinsic mathematical interest, urn models have found applications in a wide range of areas, including response-adaptive
clinical trials (\cite{Wei_Urn_clinical_trials,Asymptotic_properties_of_adaptive_designs_for_clinical_trials_with_delayed_response}), reinforcement learning (\cite{When_can_the_two-armed_bandit,A-survey-of-random-processes-with-reinforcement}), population genetics (\cite{Urn_population_genetics_peter,
Urn_population_genetics_peter_hoppe}),
social networks (\cite{social_network,urn_social_networks}), and preferential attachment models for
complex networks (\cite{preferential-attachment-and-generalized-Polya-urn,Urn_preferential_attachment}).

While classical urn models describe reinforcement within a single evolving system, many applications involve several interacting agents whose dynamics influence one another. For example, in opinion dynamics on social networks, the opinion of an individual is often influenced by the collective past opinions of its neighbors. Thus, reinforcement occurs not only through an individual's own history, but also through interactions with connected users. Interacting urn models provide a natural framework for capturing such network-based reinforcement (see, e.g., \cite{launay2011interacting,launay2012generalized,aletti2017interacting,Synchronisation-2,kaur2023interacting,urn_social_networks,Opinion_dynamics_IEEE,Social_pressure_opinion_IEEE}). More generally, interacting reinforced stochastic processes have been used to study synchronization phenomena and statistical inference for network-based reinforced systems (\cite{aletti2020statinf}). Interacting urn models with infinitely many colors have also been proposed to describe interacting innovation processes, where the emergence of novelties in one process stimulates further innovations in related processes (\cite{aletti2023interactinginnovation,statisticalinnovation}). Such models have subsequently been used to investigate technological spillovers and innovation ecosystems across different technological categories (\cite{innovationurn}). More recently, interacting urn schemes have been proposed for adaptive clinical trial designs with information sharing across different strata (\cite{aletti2024interactingurnsdesign}).

Interacting urn schemes can be broadly classified into two categories according to how the interaction is introduced. In the first, the interaction enters through the reinforcement mechanism (see \cite{kaur2023interacting,kaur2023interacting2,interacting_urns_FA-DAG}), 
which we refer to as an \emph{interacting urn scheme with reinforcement-based interaction} (see \Cref{subsec:mainmodel} for details). In the second, the interaction enters only through the drawing mechanism (see \cite{launay2011interacting,aletti2017interacting}), which we refer to as an \emph{interacting urn scheme with drawing-based interaction} (see \Cref{subsec:drawing_interaction} for details). The interaction mechanism, describing which urns influence the reinforcement or drawing mechanism of other urns, can be naturally encoded by a directed graph, where each vertex is associated with an urn and the direction of an edge indicates the direction of interaction.

A fundamental phenomenon observed in these interacting urn schemes is that, when the reinforcement matrix is irreducible, the color proportions of all urns converge almost surely to the same deterministic limit, given by the normalized Perron--Frobenius eigenvector of the reinforcement matrix (see, e.g., \cite{aletti2017interacting,kaur2023interacting,kaur2023interacting2}). A similar phenomenon is well known for classical single-urn models (see, e.g., \cite{Friedmanurn,athreya1968embedding,Randomized-urn-SA}). More generally, for a single urn, whenever the Markov chain associated with the reinforcement matrix or kernel is ergodic, the color proportions converge, and this property is not restricted to finite-color models (see \cite{bandyopadhyay2020strong,bandyopadhyay2022new}). This naturally raises the question of whether the same convergence phenomenon persists when infinitely many urns interact with one another. Addressing this question is one of the main objectives of the present article.

Our original motivation for studying infinite interacting urn schemes was to develop a simple and exactly solvable probabilistic system exhibiting \emph{self-organized criticality}, inspired by the so called Abelian sandpile model (see \cite{bak1987self,frigg2003self,pruessner2012self}). In particular, a closer examination reveals natural similarities in the dynamics of the Abelian sandpile model and interacting urn schemes with reinforcement-based interaction (see \Cref{subsec:SOC}). For the classical Abelian sandpile model, exact results are largely confined to one and two dimensions, while many properties of the model on higher-dimensional infinite lattices remain poorly understood (see \cite{DHAR1992129,dhar1993self,DHAR19994,DHAR200629,Redig2006,Jarai2018} for details).

One of the main issues in studying infinite interacting urn schemes is that, the standard tools used to analyze reinforced systems in the finite
setting, such as martingale methods and stochastic approximation (see, e.g.,
\cite{Randomized-urn-SA,kaur2023interacting,yogesh2024interacting,kaur2023interacting2}), do not extend
directly in the same generality.
Although there is a substantial literature on Banach space-valued martingales
and stochastic approximation, the available results often impose rather strong
conditions on the spectrum of the underlying reinforcement or drift operator.
These restrictions significantly limit their applicability compared with the
finite-dimensional setting.

A further difficulty is that infinitely many urns are updated simultaneously
at each time step. Consequently, there is no immediate analogue of the
Athreya--Karlin embedding, which has proved remarkably useful in the analysis
of classical urn models by connecting discrete-time urn processes with
continuous-time branching processes (see \cite{athreya1968embedding,Janson-FCLT,janson2006limit,janson2024almost}).

In this article we bypass all these difficulties by introducing a different technique which we call the \emph{representation approach}. It is described in the next subsection and in full details in \Cref{sec:representation}, is the main novelty of this work.

\subsection{The representation approach} The main goal of this article is to build an equivalent representation of interacting urn schemes that is flexible enough to work even when there are infinitely many urns. This novel idea is also a generalization of the representation approach introduced in \cite{bandyopadhyay2022new}. There, the classical urn scheme was represented by a \emph{branching Markov process} evolving in a random time environment given by the \emph{random recursive tree}.

When multiple urns interact with each other, the corresponding representation turns out to be a collection of branching Markov processes, each evolving on a tree belonging to a random forest. We shall see that this representation is good enough to study the first-order asymptotic behavior of the color proportions, even when there are infinitely many urns.

\subsection{Outline of the article}

In \Cref{subsec:mainmodel}, we describe the general model of interacting
urn schemes with reinforcement-based interaction. In
\Cref{sec:representation}, we develop the representation approach and show
how it can be used to obtain first-order asymptotics for this model, even when
the system contains infinitely many urns. In \Cref{sec:application}, we
study measure-valued generalizations and some other related models,
including interacting urn schemes with drawing-based interaction, and show how
the representation approach can be used to establish convergence in greater
generality. Finally, in \Cref{sec:proofs_Hoppe_trees}, we present the proofs
of a few intermediate lemmas used in \Cref{sec:representation}.

\subsection{Notations used in the article}
Throughout the article, all vectors are assumed to be row vectors. For a matrix $A$, we denote its transpose by $A'$. The $i$-th row vector of $A$ is denoted by $A_i$, and its $(i,j)$th entry is denoted by $A_{i,j}$. For an event $E$, $\bb{1}_E$ denotes its indicator function. We write $X\eqd Y$ if the two random variables $X$ and $Y$ have the same distribution. When $X$ and $Y$ are two discrete time stochastic processes, $X\eqd Y$ means that their finite-dimensional distributions coincide. We write $X\st Y$ if
$$
\Q(X>t)\le \Q(Y>t),\hspace{0.9cm}\forall t\in\mathbb{R},
$$
in which case we say that $X$ is stochastically dominated by $Y$. The notations $\xrightarrow{a.s.}$, $\xrightarrow{\Q}$, and $\xrightarrow{d}$ denote almost sure convergence, convergence in probability, and convergence in distribution, respectively. For a sequence of probability measures $\{\mu_n:n\ge1\}$ and a probability measure $\mu$ on a Polish space $S$, we write
$
\mu_n\implies\mu
$
if $\mu_n$ converges weakly to $\mu$, that is,
$$
\int_S f\,d\mu_n
\longrightarrow
\int_S f\,d\mu
$$
for every bounded continuous function $f:S\to\mathbb{R}$.

\section{Interacting urn schemes with reinforcement-based interaction}\label{subsec:mainmodel}
Let $\G=(V,\E)$ be a directed graph with non-empty vertex set $V$ and directed edge set $\E\subseteq V\times V$. If $(u,v)\in\E$, we say that there is a directed edge from $u$ to $v$. In this case $u$ is called an in-neighbor of $v$, and $v$ is called an out-neighbor of $u$.

For each vertex $v\in V$, let
$$
d_v^{\text{in}}
:=
\big|\{u\in V:(u,v)\in\E\}\big|,\hspace{0.9cm}
d_v^{\text{out}}
:=
\big|\{u\in V:(v,u)\in\E\}\big|
$$
and
$$
d_v:=d_v^{\text{in}}+d_v^{\text{out}}
$$
denote the in-degree, out-degree, and degree of $v$ respectively. Here $|A|$ denotes the cardinality of a set $A$.

The graph $\G$ may contain self-loops, that is, there may exist $v\in V$ such that $(v,v)\in\E$. The graph $\G$ may also be infinite. In fact, our main interest lies in the case where $\G$ is infinite. Throughout, we assume that $\G$ is locally finite, that is
$$
d_v<\infty,
\hspace{0.9cm} \forall v\in V.
$$

For $v\in V$, let
$$
N_v:=\{u\in V:(u,v)\in\E\}
$$
denote the set of in-neighbors of $v$.

Suppose that at each vertex $v\in V$ there is an urn containing balls of $k\ge2$ colors. Let
$$
\C{n}{v}
=
\left(\C{n}{v}(1),\ldots,\C{n}{v}(k)\right)
$$
denote the color-composition vector of the urn at vertex $v$ at time $n\ge0$, where $\C{n}{v}(i)$ denotes the number of balls of color $i$ present in the urn at time $n$. Further, let
$$
T_n^{(v)}
:=
\sum_{i=1}^k\C{n}{v}(i)
$$
denote the total number of balls in the urn at vertex $v$ at time $n$. The evolution of the urn configuration is defined as follows.

At each time $n\ge0$, one ball is drawn simultaneously from every urn. Let
$$
\Z{n}{v}\in\{1,\ldots,k\}
$$
denote the color drawn from the urn at vertex $v$ at time $n-1\ge0$. Conditional on
$$
\{\C{m}{u}:u\in V,\ 0\le m\le n-1\},
$$
the draws are independent and, from each urn, drawing is made uniformly at random. More precisely,
\begin{equation}\label{con}
\Q\left(
\Z{n}{v}=i_v,\ \forall v\in V_0
\,\middle|\,
\C{m}{u}:\ u\in V,\ m\le n-1
\right)
=
\prod_{v\in V_0}
\frac{\C{n-1}{v}(i_v)}{T_{n-1}^{(v)}},
\end{equation}
for every finite non-empty set $V_0\subseteq V$ and every choice of colors $i_v\in\{1,\ldots,k\}$ for all $v\in V_0$.

Reinforcement is determined by a $k\times k$ matrix
$$
R=(R_{i,j})_{1\le i,j\le k}
$$
with non-negative entries. Throughout we assume that

\begin{enumerate}[label=(2.\arabic*)]
\item\label{A1}
every urn is initially non-empty, that is,
$$
T_0^{(v)}>0,
\hspace{0.9cm} \forall v\in V;
$$

\item\label{A2}
the matrix $R$ is balanced, that is, there exists $r>0$ such that
$$
\sum_{j=1}^kR_{i,j}=r,
\hspace{0.9cm} \forall i=1,\ldots,k.
$$
\end{enumerate}

Once the colors $\{\Z{n}{v}:v\in V\}$ are sampled, all urns are updated simultaneously according to one of the following reinforcement schemes.\\[0.25cm]
\noindent\textbf{Sum-reinforcement:}

\begin{equation}\label{sum}
\C{n}{v}
=\begin{cases}
\C{n-1}{v}
+
\mathop{\sum}\limits_{u\in N_v}
R_{\Z{n}{u}},&\text{if}~~d_v^{\text{in}}>0,\\[0.25cm]
\C{n-1}{v},&\text{if}~~d_v^{\text{in}}=0.
\end{cases}
\end{equation}

\noindent\textbf{Average-reinforcement:}

\begin{equation}\label{mean}
\C{n}{v}
=\begin{cases}
\C{n-1}{v}
+
\frac{1}{d_v^{\text{in}}}
\mathop{\sum}\limits_{u\in N_v}
R_{\Z{n}{u}},&\text{if}~~d_v^{\text{in}}>0,\\[0.25cm]
\C{n-1}{v},&\text{if}~~d_v^{\text{in}}=0.
\end{cases}
\end{equation}

recall that $R_i$ denotes the $i$-th row vector of $R$. Fractional ball counts are allowed, since drawing depends only on color proportions.

The resulting process
$$
\{\C{n}{v}:v\in V,\ n\ge0\}
$$
is called an \emph{interacting urn scheme with reinforcement-based interaction, on the interaction graph $\G$, with initial configuration $\{\C{0}{v}:v\in V\}$ and reinforcement matrix $R$}. The sequence
$
\{\Z{n}{v}:v\in V,\ n\ge1\}
$
is called the color sequence of the interacting urn scheme.

Assumption \ref{A1} ensures that drawing can be performed from every urn at all times. The balancedness assumption \ref{A2} is standard in the urn literature; it implies that the total number of balls in each urn evolves linearly and deterministically. In finite systems, this property allows one to construct useful martingales and apply stochastic approximation techniques (see, e.g., \cite{Randomized-urn-SA,kaur2023interacting,yogesh2024interacting,kaur2023interacting2}).

Our ultimate goal is to study the asymptotic behavior of the color proportion vectors
$$
\frac{\C{n}{v}}{T_n^{(v)}},
$$
as well as the color-count statistics
$$
\mathcal{N}_n^{(v)}:=\frac{1}{n}
\sum_{p=1}^{n}
\left(
\bb1_{\{\Z{p}{v}=1\}},
\ldots,
\bb1_{\{\Z{p}{v}=k\}}
\right).
$$

\begin{remark}
\textnormal{For finite graph $\G$, this model with the sum-reinforcement scheme was first studied in \cite{kaur2023interacting} for the case of two colors. Assuming that the reinforcement matrix $R$ is irreducible, the authors used stochastic approximation techniques to establish almost sure convergence of the color proportion vector. Using the same framework, they also obtained rates of convergence. These results were later generalized in \cite{kaur2023interacting2}, where node-based reinforcement mechanisms were considered. A multiple-drawing version of the model was studied in \cite{yogesh2024interacting}. All of these works rely on either extracting suitable martingales or formulating an appropriate stochastic approximation scheme. However, these approaches appear to break down when $\G$ is infinite, that is, when there are infinitely many interacting urns.}
\end{remark}

\section{Representation of interacting urn schemes with reinforcement-based interaction}\label{sec:representation}
In this section, we describe in detail the representation approach for the
interacting urn scheme introduced in \Cref{subsec:mainmodel} and show how it can be
used to study the convergence of the color proportions of all urns. Our approach
is motivated by the \emph{grand representation theorem} of
\cite{bandyopadhyay2022new}. We develop a coupling between the joint urn process
and a collection of branching Markov chains indexed by a recursively growing
random time structure, which we call the \emph{interacting random recursive
forest}. We then study the asymptotic behavior of certain graph statistics
associated with this forest and establish their connection with the so called \emph{random recursive trees}.

\subsection{Random recursive trees and Hoppe trees} A \emph{uniform random recursive tree} or simply \emph{uniform recursive tree} (URT) is a random rooted tree on the vertex set $\{0,1,\ldots\}$, constructed recursively as follows. Starting with a single root vertex labeled $0$, each new vertex $j\ge1$ is connected to one of the existing vertices in $\{0,1,\ldots,j-1\}$, chosen uniformly at random and independently of the previous choices. The model was first introduced in \cite{distancenodes}; see also \cite{surveyrecursive} for a survey.

A natural generalization of the URT is the \emph{Hoppe tree}, which is associated with Hoppe's urn \cite{Hoppeurn}. Fix $t>0$. The Hoppe tree is constructed similarly to the URT, except that the root vertex $0$ is assigned weight $t$, while every other vertex is assigned weight $1$.

More precisely, let $\{D_n^t:n\ge1\}$ be a sequence of independent random variables such that $D_1^{(t)}=0$ and for $n\ge2$,
$$
\Q(D_n^t=j)
=
\begin{cases}
\displaystyle
\frac{t}{n+t-1},
&\text{if}~j=0,
\\[0.3cm]
\displaystyle
\frac{1}{n+t-1},
&\text{if}~1\le j\le n-1.
\end{cases}
$$
Let
$
\mathcal T_0^t=\{0\}
$ be the tree at time $n=0$, which contains only the root $0$. For $n\ge1$, having constructed the tree $\mathcal T_{n-1}^t$ on the vertex set $\{0,1,\ldots,n-1\}$, the tree $\mathcal T_n^t$ is obtained by attaching the vertex $n$ to the vertex $D_n^t$. The (infinite) Hoppe tree with parameter $t$ is then defined by
$$
\mathcal T^t:=\bigcup_{n\ge0}\mathcal T_n^t.
$$
The case $t=1$ corresponds to the URT.

For $n\ge 0$, let $|n|_{\mathcal T^t}$ denote the depth of the vertex $n$ in the tree $\mathcal T^t$, that is, $|n|_{\mathcal T^t}+1$ equals the number of vertices on the unique path connecting the root $0$ to $n$. The following result on the asymptotic growth of the depth $ |n|_{\mathcal T^t} $ is well known. A short proof is provided in \Cref{sec:proofs_Hoppe_trees} for completeness.

\begin{lemma}\label{lem:Hoppe_depth}
For every $t>0$,
\begin{equation}
\label{l1}
\frac{|n|_{\mathcal T^t}}{\log n}
\xrightarrow{\mathbb P}
1,\hspace{0.9cm}\text{as}~n\to\infty.
\end{equation}
\end{lemma}

Another two quantities of interest are the lowest common ancestor (LCA) and the distance between two vertices in the Hoppe tree. For vertices $m,n\ge0$, let
$
m\wedge n
$
denote the lowest common ancestor of $m$ and $n$ in $\mathcal T^t$. The graph distance between $m$ and $n$ is given by
$$
d_{\T{}{t}}(m,n)
=
|m|_{\T{}{t}}
+
|n|_{\T{}{t}}
-
2|m\wedge n|_{\T{}{t}}.
$$

It was shown in \cite{Lcaurt} that for the URT $\mathcal{T}^1$,
$$
\left|m\wedge n\right|_{\mathcal{T}^1}
\xrightarrow{d}
G,
\hspace{0.9cm}\text{as }m,n\to\infty,
$$
where $G$ is a geometric random variable with parameter $1/2$. In particular,
the family
$$
\left\{
\left|m\wedge n\right|_{\mathcal{T}^1}
:\ 1\le m<n
\right\}
$$
is tight. This tightness property is the key ingredient in establishing
convergence in probability of urn proportions in \cite{bandyopadhyay2022new}. The following lemma shows that the same property continues to hold for Hoppe trees. A proof is provided in \Cref{sec:proofs_Hoppe_trees}.

\begin{lemma}\label{lem:tightness_LCA_Hoppe_tree}
    For every $t>0$ and integer $K\ge1$,
    \begin{equation}
        \label{l2}
        \limsup_{n>m\to\infty}\Q\left(\left|m\wedge n\right|_{\mathcal{T}^t}>K\right)\le\frac{1}{2^{K-1}}.
    \end{equation}
\end{lemma}

To prove the convergence of interacting urns, we will need quantitative estimates on statistics involving distances in Hoppe trees. The following lemma, taken from \cite{bandyopadhyay2020strong}, provides the bounds that will be useful later.

\begin{lemma}[cf. Lemma 3.2 of \cite{bandyopadhyay2020strong}]
\label{lem:power_distance_Hoppe_tree}
For any $t>0$ and $\rho\in(0,1)$, there exists a constant $C>0$, depending only on $t$ and $\rho$, such that
\begin{align*}
\Ex{\sum_{n=1}^N\rho^{|n|_{\T{}{t}}}}
&\le
CN^\rho,
\\[0.15cm]
\Ex{\sum_{1\le m,n\le N}
\rho^{d_{\mathcal T^t}(m,n)}}
&\le
\begin{cases}
CN^{2\rho},
& \frac12<\rho<1,
\\[0.15cm]
CN\log N,
& \rho=\frac12,
\\[0.15cm]
CN,
& 0<\rho<\frac12.
\end{cases}
\end{align*}
\end{lemma}

Next, we formally introduce a recursive forest in which the trees evolve simultaneously according to an interactive mechanism governed by the neighboring structure of the graph $\G$. This forms the backbone of our interacting urn scheme.

\subsection{Interacting random recursive forest}\label{subsec:IRRF}
An \emph{interacting recursive forest} associated with the interaction
graph $\G$ is a forest whose vertices are of the form $(v,n)$, where
$v\in V$ and $n$ is a non-negative integer. We refer to $v$ and $n$ as
the \emph{vertex part} and \emph{time part}, respectively. The vertices
$\{(v,0):v\in V\}$ are called the \emph{roots}. For each $v\in V$ and
$n\geq1$, if $d_v^{\text{in}}>0$, the vertex $(v,n)$ is attached to a
parent vertex $(u,m)$ for some $u\in N_v$ and $0\leq m<n$. If
$d_v^{\text{in}}=0$, it is attached to the root $(v,0)$.

Observe that every path from a vertex to its root in an interacting recursive
forest is of the form
$$
(v_1,n_1)\longminus(v_2,n_2)\longminus\cdots
\longminus(v_N,n_N)\longminus(v_{N+1},0),
$$
where
$
v_{l+1}\in N_{v_l}$ and $
n_{l+1}<n_l,$ for all $1\leq l\leq N-1.
$
An interacting recursive forest involving randomness is called an \emph{interacting random recursive forest}.

Suppose $\{c_v:v\in V\}$ is a collection of positive real numbers. Let
$$
\{\D{v}{n}:v\in V,\ n\ge0\}
$$
be a family of independent random variables taking values in the set $\{(v,n):v\in V,n\ge0\}$, defined as follows.

Set
$$
\D{v}{0}=\D{v}{1}=(v,0),
$$
and for each $n\ge2$, the random variable $\D{v}{n}$ takes values in
$$
\{(v,0)\}\cup\{(u,j):u\in N_v,\ 1\le j\le n-1\}
$$
such that if $d_v^{\text{in}}=0$ then $\D{v}{n}=(v,0)$ and if $d_v^{\text{in}}>0$ then the probability distribution is given by
\begin{align}
&\Q\left(\D{v}{n}=\gamma\right)\notag\\[0.25cm]
=&\begin{cases}
\displaystyle
\frac{c_v}{c_v+(n-1)d_v^{\text{in}}},
& \text{if }\gamma=(v,0),
\\[0.3cm]
\displaystyle
\frac{1}{c_v+(n-1)d_v^{\text{in}}},
& \text{if }\gamma=(u,j),\ u\in N_v,\ 1\le j\le n-1.\label{eq:Dvn_distribution_type1}
\end{cases}
\end{align}

Using these variables, we recursively construct a forest
$$
\Fo{}=\bigcup_{n\ge0}\Fo{n}
$$
with vertex set
$
\{(v,n):v\in V,\ n\ge0\}
$ as follows.

The initial forest $\Fo{0}$ consists of only the root vertices
$
\{(v,0):v\in V\}.
$
For each $n\ge1$, having constructed $\Fo{n-1}$, we obtain $\Fo{n}$ by adding the $n$'th generation vertices
$
\{(v,n):v\in V\}
$
connecting each $(v,n)$ to its parent $\D{v}{n}$.

Each connected component of $\Fo{}$ is a rooted tree whose root is of the form $(v,0)$ for some $v\in V$. By construction it is immediate that any realization of $\mathbb F$ is an interacting recursive forest associated with $\G$. We call $\Fo{}$ the \emph{interacting random recursive forest of first kind} associated with the interaction graph $\G$ and root weights $\{c_v:v\in V\}$, and write
$$
\Fo{}
\sim
\operatorname{IRRF}_1\left(\G,\{c_v:v\in V\}\right).
$$

We next define a second forest
$$
\Fo{}'=\bigcup_{n\ge0}\Fo{n}',
$$
constructed in the same way, except that the parent variables satisfy
$$
\D{v}{0}=\D{v}{1}=(v,0)
$$
and if $n\ge2$ and $d_v^{\text{in}}>0$ then
\begin{align}
&\Q\left(\D{v}{n}=\gamma\right)\notag\\[0.25cm]
=&
\begin{cases}
\displaystyle
\frac{c_v}{c_v+n-1},
& \text{if }\gamma=(v,0),
\\[0.3cm]
\displaystyle
\frac{1}{(c_v+n-1)d_v^{\text{in}}},
& \text{if }\gamma=(u,j),\ u\in N_v,\ 1\le j\le n-1.\label{eq:Dvn_distribution_type2}
\end{cases}
\end{align}
and if $d_v^{\text{in}}=0$ then $\D{v}{n}=(v,0)$ for all $n\ge2$.

The resulting forest is called the \emph{interacting random recursive forest of second kind}, and we write
$$
\Fo{}'
\sim
\operatorname{IRRF}_2\left(\G,\{c_v:v\in V\}\right).
$$

As an example, consider the directed line
$$
V=\mathbb Z,
\hspace{0.9cm}
\E=\{(v,v+1):v\in\mathbb Z\}.
$$
In this case,
$$
N_v=\{v-1\},
\hspace{0.9cm} \forall v\in\mathbb Z,
$$
and therefore the vertex $(v,n)$ either attaches to the root $(v,0)$ or to one of the previously created vertices
$
\{(v-1,j):1\le j\le n-1\}.
$

\begin{figure}[h]
\centering

\begin{tikzpicture}[scale=0.82]

\tikzset{
edge/.style={thick},
every node/.style={font=\small}
}

\node at (-7.2,0) {$\cdots$};

\node (m3) at (-6,0) {$(-3,0)$};
\node (m2) at (-3,0) {$(-2,0)$};
\node (m1) at (0,0) {$(-1,0)$};
\node (z0) at (3,0) {$(0,0)$};
\node (p1) at (6,0) {$(1,0)$};

\node at (7,0) {$\cdots$};

\draw[->,edge] (-5.2,0) -- (-3.8,0);
\draw[->,edge] (-2.2,0) -- (-0.8,0);
\draw[->,edge] (0.8,0) -- (2.2,0);
\draw[->,edge] (3.8,0) -- (5.2,0);

\node (a1) at (-6,1.8) {$(-3,1)$};
\node (a2) at (-7,3.4) {$(-2,3)$};
\node (a3) at (-5,3.4) {$(-2,2)$};

\draw[edge] (a1) -- (m3);
\draw[edge] (a2) -- (a1);
\draw[edge] (a3) -- (a1);

\node at (-7,4.7) {$\vdots$};
\node at (-5,4.7) {$\vdots$};

\node (b1) at (-3,1.8) {$(-2,1)$};
\node (b2) at (-3,3.4) {$(-1,3)$};

\draw[edge] (b1) -- (m2);
\draw[edge] (b2) -- (b1);

\node at (-3,4.7) {$\vdots$};

\node (c1) at (0,1.8) {$(-1,1)$};
\node (c2) at (1.2,3.4) {$(0,2)$};
\node (c3) at (1.5,1.8) {$(-1,2)$};

\draw[edge] (c1) -- (m1);
\draw[edge] (c2) -- (c1);
\draw[edge] (c3) -- (m1);

\node at (1.2,4.7) {$\vdots$};
\node at (1.5,2.8) {$\vdots$};

\node (d1) at (3,1.8) {$(0,1)$};
\node (d3) at (4.5,1.8) {$(0,3)$};

\draw[edge] (d1) -- (z0);
\draw[edge] (d3) -- (z0);

\node at (3,3.1) {$\vdots$};
\node at (4.5,3.1) {$\vdots$};

\end{tikzpicture}

\caption{
A realization of the interacting random recursive forest on the directed line:
$
V=\mathbb Z
$,
$
\E=\{(v,v+1):v\in\mathbb Z\}
$.
}
\label{fig:irrf_line}
\end{figure}
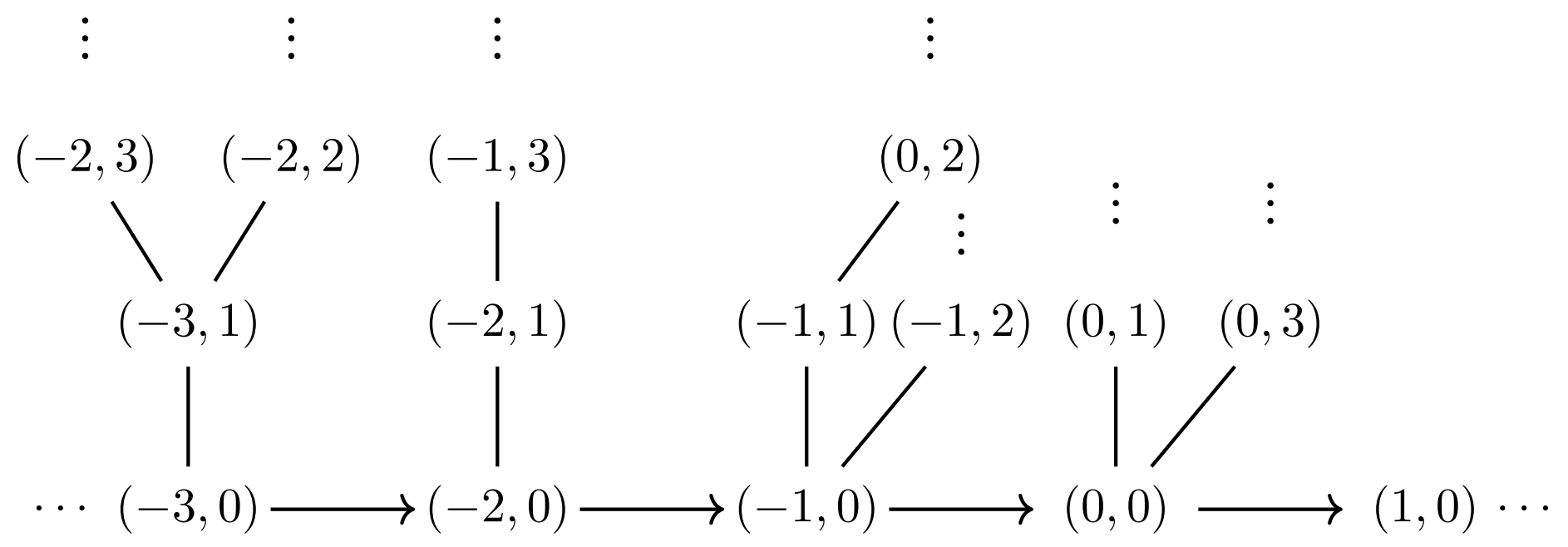

\Cref{fig:irrf_line} shows a realization of several trees in the interacting random recursive forest on the directed line. The horizontal arrows represent the interaction graph, while the vertices above each root illustrate the corresponding trees generated by the attachment mechanism. In the realization shown in \Cref{fig:irrf_line}, we have for the first generation vertices
\begin{align*}
&\ldots\ldots~\D{-3}{1}=(-3,0),\hspace{0.9cm}
\D{-2}{1}=(-2,0),\hspace{0.9cm}\D{-1}{1}=(-1,0),\\[0.1cm]
&\hspace{0.9cm}\hspace{0.9cm}\hspace{0.9cm}\hspace{0.9cm}\D{0}{1}=(0,0),~\ldots\ldots
\end{align*}
for the second generation vertices
$$
\ldots~\D{-2}{2}=(-3,1),
\hspace{0.9cm}
\D{-1}{2}=(-1,0),
\hspace{0.9cm}
\D{0}{2}=(-1,1),~\ldots
$$
and for the third generation vertices
$$
\ldots~\D{-2}{3}=(-3,1),
\hspace{0.9cm}
\D{-1}{3}=(-2,1),
\hspace{0.9cm}
\D{0}{3}=(0,0),~\ldots
$$

 We next introduce several graph-theoretic statistics associated with an interacting recursive forest $\mathbb F$. Let $\mathcal R((v,n))$ denote the root of the tree containing $(v,n)$. The depth of a vertex $(v,n)$, denoted by $|(v,n)|_{\mathbb F}$, is defined by the graph distance between $(v,n)$ and $\R{(v,n)}$ in $\mathbb F$. 
 
 If
$
\mathcal R((u,m))
=
\mathcal R((v,n)),
$
that is $(u,m)$ and $(v,n)$ are in same tree of $\mathbb F$, then we write
$
(u,m)\wedge(v,n)
$
for the lowest common ancestor of $(u,m)$ and $(v,n)$. To unify notation, we set
$$
|(u,m)\wedge(v,n)|_{\mathbb F}:=0
$$
whenever
$
\mathcal R((u,m))
\neq
\mathcal R((v,n)).
$

The distance between two vertices $(u,m)$ and $(v,n)$ in $\mathbb F$ is defined by
$$
d_{\mathbb F}((u,m),(v,n))
=
|(u,m)|_{\mathbb F}
+
|(v,n)|_{\mathbb F}
-
2|(u,m)\wedge(v,n)|_{\mathbb F}.
$$
When the vertices belong to different trees, this reduces to
$$
d_{\mathbb F}((u,m),(v,n))
=
|(u,m)|_{\mathbb F}
+
|(v,n)|_{\mathbb F}.
$$
It is straightforward to verify that $d_{\mathbb F}$ defines a metric on
$
\{(v,n):v\in V,\ n\ge0\}.
$ 

\subsection{Branching Markov process on interacting random recursive forest}\label{subsec:branchingmarkovprocess}

Let $u\in[0,1)$ and let
$
p=(p_1,\ldots,p_k)
$
be a probability vector, that is,
$$
p_i\ge0
\hspace{0.9cm}\text{and}\hspace{0.9cm}
\sum_{i=1}^k p_i=1.
$$
Setting $p_0=0$, define the function
$$
f(u,p)
=
\sum_{j=1}^k
j\,
\mathbbm{1}_{\left\{
\sum_{i=0}^{j-1}p_i
\le
u
<
\sum_{i=0}^{j}p_i
\right\}}.
$$
Thus, if $U$ is uniformly distributed on $[0,1)$, then $f(U,p)$ is a random variable taking values in $\{1,\ldots,k\}$ with distribution $p$.

Let
$$
\Fo{}=\bigcup_{n\ge0}\Fo{n}
$$
be an interacting random recursive forest (either of first or second kind) associated with the interaction graph $\G$ and root weights
$
\{r^{-1}T_0^{(v)}:v\in V\}.
$
Further, let
$
\{U_{(v,n)}:v\in V,\ n\ge1\}
$
be a collection of i.i.d.\ uniform random variables on $[0,1)$, independent of the forest $\Fo{}$.

We now define a process
$
\{\W{v}{n}:v\in V,\ n\ge1\}
$
on the vertices of $\Fo{}$, where $\W{v}{n}\in\{1,\ldots,k\}
$
denotes the color assigned to the vertex $(v,n)$.

For $v\in V$, let
$$
\W{v}{1}
:=
f\left(
U_{(v,1)},
\C{0}{v}/T_0^{(v)}
\right)
$$

and recursively, for $n\ge2$, define
\be\label{eq:def_Wvn}
\W{v}{n}
=
\begin{cases}
\displaystyle
f\left(
U_{(v,n)},
\C{0}{v}/T_0^{(v)}
\right),
& \text{if }\D{v}{n}=(v,0),
\\[0.25cm]
\displaystyle
f\left(
U_{(v,n)},
r^{-1}R_{W_{\D{v}{n}}}
\right),
& \text{if}~\D{v}{n}\neq(v,0).
\end{cases}
\ee

Thus, conditioned on the forest $\Fo{}$, the color of a vertex is determined by the color of its parent through the transition matrix $r^{-1}R$. Vertices attached directly to a root instead receive colors according to the initial distribution $\C{0}{v}/T_0^{(v)}.$ Hence along any ancestral path $P$ of the forest $\Fo{}$, the sequence
$
\{W_a:a\in P\}
$
forms a Markov chain on the state space $\{1,\ldots,k\}$ with transition matrix $r^{-1}R$. Consequently, the process
$
\{\W{v}{n}:v\in V,\ n\ge1\}
$
may be viewed as a branching Markov process on the interacting random recursive forest $\Fo{}$.
\subsection{The main representation theorem}

The following theorem provides a representation of interacting urn schemes with reinforcement-based interaction as branching Markov process on interacting random recursive forests.

\begin{theorem}\label{thm:representation_theorem}
Suppose
$
\{\C{n}{v}:v\in V,\ n\ge0\}
$
is an interacting urn scheme defined in \Cref{subsec:mainmodel} with interaction graph $\G=(V,\E)$, initial configuration
$
\{\C{0}{v}:v\in V\},
$
and reinforcement matrix $R$ such that assumptions \ref{A1} and \ref{A2} are satisfied. 

Then the associated color sequence
$
\{\Z{n}{v}:v\in V,\ n\ge1\}
$
satisfies
\begin{equation}
\label{t1}
\left\{\Z{n}{v}:v\in V,\ n\ge1\right\}
\stackrel{d}{=}
\left\{\W{v}{n}:v\in V,\ n\ge1\right\},
\end{equation}
where
$
\{\W{v}{n}:v\in V,\ n\ge1\}
$
is the branching Markov process on an interacting random recursive forest $\Fo{}$ defined in \Cref{subsec:branchingmarkovprocess}. The forest is of the first or second kind, depending on whether the interacting urn scheme uses a sum-reinforcement or an average-reinforcement scheme.
\end{theorem}

\begin{proof}
Since all row sums of $R$ is same, it can be checked that $T_n^{(v)}=T_0^{(v)}+nrd_v^{\text{in} }$ for the sum-reinforcement scheme and $T_n^{(v)}=T_0^{(v)}+nr$ for the average-reinforcement scheme.

For $v\in V$, define
$$
H_{(v,0)}:=\C{0}{v},
$$

and for $n\ge1$, if $d_v^{\text{in}}=0$ then set $H_{(v,n)}=H_{(v,0)}$ and if $d_v^{\text{in}}>0$ then 
$$
H_{(v,n)}
:=
H_{(v,0)}
+
\begin{cases}
\displaystyle
\sum_{j=1}^n
\sum_{u\in N_v}
R_{\W{u}{j}},
&
\text{if }\Fo{}\sim
\operatorname{IRRF}_1
\left(
\G,{\{r^{-1}T_0^{(v)}\}}
\right),
\\[0.6cm]
\displaystyle
\frac{1}{d_v^{\text{in}}}
\sum_{j=1}^n
\sum_{u\in N_v}
R_{\W{u}{j}},
&
\text{if }\Fo{}\sim
\operatorname{IRRF}_2
\left(
\G,
\{r^{-1}T_0^{(v)}\}
\right).
\end{cases}
$$

We now prove by induction that
\begin{equation}
\label{eq:induction_hypothesis}
\left\{
(\W{v}{j},H_{(v,j-1)}):
v\in V,\ j\le n
\right\}
\stackrel{d}{=}
\left\{
(\Z{j}{v},\C{j-1}{v}):
v\in V,\ j\le n
\right\},
\end{equation}
for all $n\ge1$.

A quick observation is that if \eqref{eq:induction_hypothesis} holds for some $n$ then necessarily
$$
\{H_{(v,n)}:v\in V\}\eqd\{\C{n}{v}:v\in V\}.
$$
For $n=1$, by definition,
$
\W{v}{1}
$
has distribution
$
\C{0}{v}/T_0^{(v)},
$
which coincides with the distribution of
$
\Z{1}{v}.
$
Since
$
H_{(v,0)}=\C{0}{v},
$
\eqref{eq:induction_hypothesis} follows for $n=1$.

Now assume that \eqref{eq:induction_hypothesis} holds for some $n=m-1\ge1$. Define
$$
\mathcal G_{m-1}
:=
\sigma(\W{v}{j}:v\in V,\ j\le m-1).
$$
From the definition \eqref{eq:def_Wvn}, for each $v\in V$ and $1\le i\le k$,
\begin{align*}
&\Q\left(
\W{v}{m}=i
\,\middle|\,
\Fo m,\mathcal G_{m-1}
\right)\\[0.25cm]
=&\frac{H_{(v,0)}(i)}{T_0^{(v)}}
\,
\bb 1_{\{\D{v}{m}=(v,0)\}}
+
r^{-1}R_{W_{\D{v}{m}},i}
\,
\bb 1_{\{\D{v}{m}\ne(v,0)\}}.
\end{align*}
 Notice that, conditioned on the forest $\Fo m$ and $\mathcal G_{m-1}$, the random variables
$
\{\W{v}{m}:v\in V\}
$
are independent. Hence, for any finite non-empty subset $V_0\subseteq V$ and any collection
$
\{i^{(v)}:v\in V_0\}
\subseteq
\{1,\ldots,k\}
$
of colors, we have
\begin{align}
&
\Q\Big(
\W{v}{m}=i^{(v)},
\,v\in V_0
\,\Big|\,
\Fo m,\mathcal G_{m-1}
\Big)\notag
\\[0.25cm]
=&
\prod_{v\in V_0}
\left[
\frac{H_{(v,0)}(i^{(v)})}{T_0^{(v)}}
\,
\bb 1_{\{\D{v}{m}=(v,0)\}}
+
r^{-1}R_{W_{\D{v}{m}},i^{(v)}}
\,
\bb 1_{\{\D{v}{m}\ne(v,0)\}}
\right]
.\label{eq:conditional_expression}
\end{align}
Using the distributions of $\D{v}{m}$ given in
\eqref{eq:Dvn_distribution_type1}
and
\eqref{eq:Dvn_distribution_type2},
a direct computation shows that, the conditional expectation
\begin{align*}
&\Ex{
\frac{H_{(v,0)}(j)}{T_0^{(v)}}
\,
\bb 1_{\{\D{v}{m}=(v,0)\}}
+
r^{-1}R_{W_{\D{v}{m}},j}
\,
\bb 1_{\{\D{v}{m}\ne(v,0)\}}~\Bigg|~\mathcal{G}_{m-1}}\\[0.25cm]
=&\frac{H_{(v,m-1)}(j)}{\sum_{l=1}^kH_{(v,m-1)}(l)}.
\end{align*}
for both sum- and average-reinforcement schemes.

Taking expectation with respect to the forest $\Fo m$ in \eqref{eq:conditional_expression}, and using the independence of
$
\{\D{v}{m}:v\in V\}
$
from $\mathcal G_{m-1}$, we obtain
$$
\Q\Big(
\W{v}{m}=i^{(v)},
\,v\in V_0
\,\Big|\,
\mathcal G_{m-1}
\Big)
=
\prod_{v\in V_0}
\frac{H_{(v,m-1)}(i^{(v)})}{\sum_{l=1}^kH_{(v,m-1)}(l)}
.
$$
which coincides with the conditional color distribution in the interacting urn dynamics \eqref{con}.

Therefore,
$
\{\W{v}{j}:v\in V,\ j\le m\}
$
and
$
\{\Z{j}{v}:v\in V,\ j\le m\}
$
have the same distribution. Further since $$\{H_{(v,m-1)}:v\in V\}\eqd\{\C{m-1}{v}:v\in V\},$$
this immediately yields
$$
\left\{
(\W{v}{j},H_{(v,j-1)}):
v\in V,\ j\le m
\right\}
\stackrel{d}{=}
\big\{
(\Z{j}{v},\C{j-1}{v}):
v\in V,\ j\le m
\big\}.
$$
This completes the induction and hence the proof.
\end{proof}

\begin{remark}
\textnormal{We mention that \Cref{thm:representation_theorem} is true for any balanced reinforcement matrix $R$. No irreducibility assumption is needed.}
\end{remark}

\begin{remark}
\textnormal{\Cref{thm:representation_theorem} extends Theorem 2.4 of \cite{bandyopadhyay2022new}, which was established for a single urn with colors taking values in a Polish space. The single urn setting is the case when $\G$ consists of a single vertex with a self-loop, in which case the underlying random forest is simply a Hoppe tree and the color sequence forms a branching Markov process on that tree. For a general graph $\G$, especially when $\G$ is infinite, the underlying forest may contain infinitely many trees, with a branching Markov process evolving on each tree.}
\end{remark}

\subsection{Properties of certain statistics in interacting random recursive forests}\label{subsec:properties_of_IRRF_statistics}
\newcommand{\uF}{\underline{\mathbb{F}}}
\newcommand{\oF}{\overline{\mathbb{F}}}
\newcommand{\bF}{\mathbb{F}}
\newcommand{\brackat}[1]{\left(#1\right)}
\newcommand{\sh}[1]{{#1}^{(-1)}}
In this section we study the asymptotic behavior of several important graph statistics associated with interacting random recursive forests. For a general interaction graph $\G$ and arbitrary root weights $\{c_v:v\in V\}$, analyzing these statistics directly turns out to be difficult. Instead, we first introduce a simpler forest where the asymptotic behavior can be studied more explicitly through a coupling with Hoppe trees. We then compare the statistics of the original forest with those of the simpler forest to derive the desired asymptotic properties. The construction of this simpler forest is given below.

Let $\{c_v:v\in V\}$ be a set of positive root weights as before. Let
$$
\{I_{(v,n)},J_{(v,n)}:v\in V,\ n\ge0\}
$$
be a collection of independent random variables defined as follows. 

For each $v\in V$ and $n\ge0$, the random variable $I_{(v,n)}$ satisfies
\begin{equation}
\label{eq:def_Ivn}
\Q(I_{(v,n)}=u)=
\begin{cases}
\displaystyle
\frac{1}{d_v^{\text{in}}}\bb{1}_{\{u\in N_v\}},
& \text{if } d_v^{\text{in}}>0,
\\[0.25cm]
\bb{1}_{\{u=v\}},
& \text{if } d_v^{\text{in}}=0.
\end{cases}
\end{equation}
Also, let
$$
J_{(v,0)}=J_{(v,1)}=0,
$$
and for $n\ge2$,
\begin{equation}
\label{eq:def_J}
\Q\left(J_{(v,n)}=j\right)
=
\begin{cases}
\displaystyle
\frac{c_v}{c_v+n-1},
& \text{if } j=0,
\\[0.3cm]
\displaystyle
\frac{1}{c_v+n-1},
& \text{if } 1\le j\le n-1.
\end{cases}
\end{equation}

Define the parent variables
$$
D_{(v,n)}=\begin{cases}
(I_{(v,n)},J_{(v,n)}),&\text{if}~d_v^{\text{in}}>0,\\[0.25cm]
(v,0),&\text{if}~d_v^{in}=0.
\end{cases}
$$
We construct the forest 
$$
\bF:=\bigcup_{n\ge0}\bF_n
$$
recursively just like in \Cref{subsec:IRRF}, as follows. The initial forest is
$$
\bF_0:=\{(v,0):v\in V\},
$$
which contains only the root vertices. For $n\ge1$, having constructed $\bF_{n-1}$, we obtain $\bF_n$ by attaching the vertex $(v,n)$ to $D_{(v,n)}$ for every $v\in V$. 

The key property of the forest $\bF$ is that, when $d_v^{\text{in}}>0$ for all $v$, the components of the parent vertices $D_{(v,n)}$ are independent.
 
 By construction $\bF$ is also an interacting random recursive forest. We call it \emph{interacting random recursive forest of third kind} associated with the interaction graph $\G$ and root weights $\{c_v:v\in V\}$, and write
$$
\bF\sim\text{IRRF}_3\big(\G,\{c_v:v\in V\}\big).
$$
In particular if $c_v=c$ for all $v\in V$, then we simply write
$$
\bF^c\sim\text{IRRF}_3\left(\G,c\right).
$$
The branching Markov process $\{\W{v}{n}:v\in V,\ n\geq0\}$ on $\bF$
can be defined in exactly the same way as in
\Cref{subsec:branchingmarkovprocess}.

We now introduce a class of events that depend only on non-root attachments in any interacting random recursive forest associated with $\G$.

\begin{definition}\label{def1}
An event $A$ defined on an interacting random recursive forest associated with the graph $\G$ is called a \emph{non-root configuration} if there exist finitely many vertices
$$
(v_1,n_1),\ldots,(v_m,n_m),
$$
with $n_i\ge2$ for every $1\le i\le m$, and vertices
$$
(u_1,p_1),\ldots,(u_m,p_m),
$$
such that
$$
u_i\in N_{v_i},
\hspace{0.9cm}
1\le p_i\le n_i-1,
$$
for every $1\le i\le m$, and
$$
A
=
\bigcap_{i=1}^m
\left\{
D_{(v_i,n_i)}
=
(u_i,p_i)
\right\}.
$$
\end{definition}

Thus a non-root configuration specifies finitely many attachment relations that do not involve any root vertex in the corresponding interacting random recursive forest.

For a non-root configuration $A$, we define its \emph{1-shift} by
$$
A^{(-1)}
:=
\bigcap_{i=1}^m
\left\{
D_{(v_i,n_i-1)}
=
(u_i,p_i-1)
\right\}.
$$

Observe that $A^{(-1)}$ need not itself be a non-root configuration, since some of the indices $p_i-1$ may equal $0$.

\begin{definition}
An event $A$ defined on an interacting random recursive forest associated with the graph $\G$ is called \emph{admissible} if it can be written as a countable disjoint union of non-root configurations.
\end{definition}

If
$$
A=\bigcup_{j=1}^\infty A_j,
$$
where $\{A_j\}_{j\ge1}$ are mutually disjoint non-root configurations, we define the 1-shift of $A$ by
$$
A^{(-1)}
:=
\bigcup_{j=1}^\infty A_j^{(-1)}.
$$

We write $\Q_{\bF}(A)$ for the probability of the event $A$ under the interacting random recursive forest $\bF$.

\begin{lemma}\label{lem:admissible_comparison}
Suppose $\bF$ is an interacting random recursive forest associated with
$\G$, of the first, second, or third kind, with positive root weights
$\{c_v:v\in V\}$ satisfying
$$
c:=\sup_{v\in V}c_v<\infty.
$$

Then for any admissible event $A$,
$$
\Q_{\bF^c}(A)
\le\Q_{\bF}(A)\le
\Q_{\bF^1}(\sh{A}).
$$
\end{lemma}

\begin{proof}
By countable additivity, it suffices to prove the result only for non-root configurations. Let
$$
A
=
\bigcap_{i=1}^m
\left\{
D_{(v_i,n_i)}
=
(u_i,p_i)
\right\},
$$
where $u_i,v_i,n_i$ and $p_i$ are as in \Cref{def1}. If $d_{v_i}^{\text{in}}=0$ for some $v_i$ then $\Q_{\mathbb F}(A)=0$ and the lemma holds trivially, so without loss we assume $d_{v_i}^{\text{in}}>0$ for all $i$.

Since the parent variables
$
\{D_{(v_i,n_i)}:1\le i\le m\}
$
are independent, it follows from
\eqref{eq:Dvn_distribution_type1}
and
\eqref{eq:Dvn_distribution_type2}
that
$$
\Q_{\bF}(A)
=
\begin{cases}
\displaystyle
\prod_{i=1}^m
\frac{1}{
c_{v_i}
+
(n_i-1)d_{v_i}^{\text{in}}
},
& \text{If}~\bF~\text{is of first kind},
\\[0.5cm]
\displaystyle
\prod_{i=1}^m
\frac{1}{
(c_{v_i}+n_i-1)
d_{v_i}^{\text{in}}
},
& \text{If}~\bF~\text{is of second or third kind}.
\end{cases}
$$

Since
$
c_{v_i}\le c,
$
we have
$$
\frac{1}{
(c+n_i-1)d_{v_i}^{\text{in}}
}
\le
\frac{1}{
c_{v_i}
+
(n_i-1)d_{v_i}^{\text{in}}
}
\le
\frac{1}{
(n_i-1)d_{v_i}^{\text{in}}
}\quad\text{and}
$$
$$
\frac{1}{
(c+n_i-1)d_{v_i}^{\text{in}}
}
\le
\frac{1}{
(c_{v_i}
+
n_i-1)d_{v_i}^{\text{in}}
}
\le
\frac{1}{
(n_i-1)d_{v_i}^{\text{in}}
}.
$$

Observe that
$$
\frac{1}{
(c+n_i-1)d_{v_i}^{\text{in}}
}
$$
is precisely the probability that
$
D_{(v_i,n_i)}=(u_i,p_i)
$
in the forest $\bF^c$, while
$$
\frac{1}{
(n_i-1)d_{v_i}^{\text{in}}
}
$$
is the probability that
$
D_{(v_i,n_i-1)}=(u_i,p_i-1)
$
in the forest $\bF^1$.

Multiplying over $1\le i\le m$ yields
$$
\Q_{\bF^c}(A)
\le
\Q_{\bF}(A)
\le
\Q_{\bF^1}(\sh{A}),
$$
which completes the proof.
\end{proof}

The following lemma compares the behavior of several graph statistics in different interacting random recursive forests associated with $\G$.

\begin{lemma}\label{lem:comparison_with_IRRF2}
Suppose $\bF$ is an interacting random recursive forest associated with
$\G$, of the first, second, or third kind, with positive root weights
$\{c_v:v\in V\}$ satisfying
$$
c:=\sup_{v\in V}c_v<\infty.
$$
Then, for every $u,v\in V$ and $m,n\ge1$, the following relations hold:
\begin{align}
|(v,n)|_{\bF^c}
&\st
|(v,n)|_{\bF}
\st
|(v,n-1)|_{\bF^1}+1,
\label{eq:depth_compare_general}
\\
|(u,m)\wedge(v,n)|_{\bF}
&\st
|(u,m-1)\wedge(v,n-1)|_{\bF^1}+1,
\label{eq:lca_compare_general}
\\
\rho^{d_{\bF}((u,m),(v,n))}
\bb{1}_{\{|(u,m)\wedge(v,n)|_{\bF}>0\}}
&\st
\rho^{d_{\bF^1}((u,m-1),(v,n-1))}
\bb{1}_{\{\R{(u,m)}=\R{(v,n)}\}},
\label{eq:distance_compare_general}
\end{align}
where $\rho\in(0,1)$ is a constant.
\end{lemma}

\begin{proof}
We first prove \eqref{eq:depth_compare_general} and
\eqref{eq:lca_compare_general}. If either $m=1$ or $n=1$, then both
statements are immediate. Hence, throughout this part of the proof we assume
$m,n\ge2$.

Fix an integer $K\ge0$. For notational simplicity, we suppress the
subscripts $\bF$, $\bF^c$, and $\bF^1$ from the depth and distance notation. Let
$$
(u_0,m_0)=(u,m),
\hspace{0.9cm}
(v_0,n_0)=(v,n).
$$

To prove \eqref{eq:depth_compare_general}, it is enough to show
\begin{equation}
\label{eq:xgsb}
\Q_{\bF^c}(|(v,n)|>K)
\le
\Q_{\bF}(|(v,n)|>K)
\le
\Q_{\bF^1}(|(v,n-1)|+1>K).
\end{equation}
The inequalities are trivial when $K=0$, since each probability is equal
to one. Thus, we may assume $K\ge1$.

Consider the event
$$
A:=\{|(v,n)|>K\}.
$$
The event $A$ occurs if and only if

\begin{enumerate}
\item
$d_v^{\text{in}}>0$;

\item
there exist vertices
$
v_1,\ldots,v_K,
$
each having positive in-degree, such that
$$
v_1\in N_{v_0},
\quad
v_2\in N_{v_1},
\quad
\ldots,
\quad
v_K\in N_{v_{K-1}};
$$

\item
there exist integers
$
n_0>n_1>\cdots>n_K\ge1,
$
such that
$$
D_{(v_{i-1},n_{i-1})}
=
(v_i,n_i),
\hspace{0.9cm}
i=1,\ldots,K.
$$
\end{enumerate}

This is precisely the non-root configuration shown in \Cref{fig:depth_case}.

\begin{figure}[h]
\centering
\begin{tikzpicture}[scale=1.1]
\tikzset{
edge/.style={thick},
every node/.style={font=\small}
}

\node (u1) at (0,1) {$(v_0,n_0)$};
\node (u2) at (2,1) {$(v_1,n_1)$};
\node (up) at (6,1) {$(v_{K},n_{K})$};

\draw (u1)--(u2);
\draw (u2)--(3.5,1);
\node at (4,1) {$\cdots$};
\draw (4.5,1)--(up);

\end{tikzpicture}
\caption{A non-root configuration corresponding to the event
$ |(v,n)|>K $.}
\label{fig:depth_case}
\end{figure}

Since the initial vertex $(v,n)$ and the integer $K$ are fixed, all
such configurations are mutually disjoint. Therefore, $A$ is admissible.
Applying \Cref{lem:admissible_comparison}, together with the identity
$$
A^{(-1)}
=
\{|(v,n-1)|\ge K\}
=
\{|(v,n-1)|+1>K\},
$$
immediately yields \eqref{eq:xgsb}, thereby proving
\eqref{eq:depth_compare_general}.

Next we prove \eqref{eq:lca_compare_general}. If
$(u,m)=(v,n)$, then the result follows immediately from
\eqref{eq:depth_compare_general}. Hence we assume
$
(u,m)\neq(v,n)
$ and $m,n\ge2$.

Consider the event
$$
B:=\{|(u,m)\wedge(v,n)|>K\}.
$$
We shall show that $B$ is admissible. The event $B$ occurs if and only if
one of the following two cases holds.
\vspace{0.25cm}

\noindent
\textit{Case 1.} $(u,m)$ does not lie on the ancestral line of $(v,n)$.

In this case, the forest must contain a non-root configuration shown in \Cref{fig:LCA_case1}, where $i_K>0$.

\begin{figure}[h]
\centering
\begin{tikzpicture}[scale=0.92]
\tikzset{
edge/.style={thick},
every node/.style={font=\small}
}

\node (u1) at (0,1) {$(u_0,m_0)$};
\node (u2) at (2,1) {$(u_1,m_1)$};
\node (up) at (6,1) {$(u_p,m_p)$};

\node (v1) at (0,0) {$(v_0,n_0)$};
\node (v2) at (2,0) {$(v_1,n_1)$};
\node (vq) at (6,0) {$(v_q,n_q)$};

\node (w1) at (8,1) {$(w_0,i_0)$};
\node (wk) at (12,1) {$(w_K,i_K)$};

\draw (u1)--(u2);
\draw (u2)--(3.5,1);
\node at (4,1) {$\cdots$};
\draw (4.5,1)--(up);

\draw (v1)--(v2);
\draw (v2)--(3.5,0);
\node at (4,0) {$\cdots$};
\draw (4.5,0)--(vq);

\draw (up)--(w1);
\draw (vq)--(w1);

\draw (w1)--(9.3,1);
\node at (10,1) {$\cdots$};
\draw (10.7,1)--(wk);

\end{tikzpicture}
\caption{A non-root configuration corresponding to the event
$ |(u,m)\wedge(v,n)|>K $ when $(u,m)$ is not an ancestor of $(v,n)$.}
\label{fig:LCA_case1}
\end{figure}

\vspace{0.25cm}

\noindent
\textit{Case 2.} $(u,m)$ lies on the ancestral line of $(v,n)$.

In this case, the forest must contain a non-root configuration shown in \Cref{fig:LCA_case2}, where again $i_K>0$.

\begin{figure}[h]
\centering
\begin{tikzpicture}[scale=0.9]
\tikzset{
edge/.style={thick},
every node/.style={font=\small}
}

\node (u1) at (0,1) {$(v_0,n_0)$};
\node (u2) at (2,1) {$(v_1,n_1)$};
\node (up) at (6,1) {$(u_0,m_0)$};

\node (w1) at (8,1) {$(w_0,i_0)$};
\node (wk) at (12,1) {$(w_K,i_K)$};

\draw (u1)--(u2);
\draw (u2)--(3.5,1);
\node at (4,1) {$\cdots$};
\draw (4.5,1)--(up);

\draw (up)--(w1);

\draw (w1)--(9.3,1);
\node at (10,1) {$\cdots$};
\draw (10.7,1)--(wk);

\end{tikzpicture}
\caption{A non-root configuration corresponding to the event
$ |(u,m)\wedge(v,n)|>K $ when $(u,m)$ is an ancestor of $(v,n)$.}
\label{fig:LCA_case2}
\end{figure}

In both figures, the corresponding event is a non-root configuration. Since the initial vertices $(u_0,m_0)$, $(v_0,n_0)$, and the integer $K$ are fixed, all such configurations are mutually disjoint. Therefore, $B$ is admissible.

Applying \Cref{lem:admissible_comparison} and observing that
$$
B^{(-1)}
=
\{|(u,m-1)\wedge(v,n-1)|\ge K\}
=
\{|(u,m-1)\wedge(v,n-1)|+1>K\},
$$
we obtain
$$
\Q_{\bF}\left(|(u,m)\wedge(v,n)|>K\right)
\le
\Q_{\bF^1}\left(|(u,m-1)\wedge(v,n-1)|+1>K\right),
$$
which proves \eqref{eq:lca_compare_general}.

Finally, we now prove \eqref{eq:distance_compare_general}. If
$(u,m)=(v,n)$, then the result is immediate. Hence, throughout the remainder of the proof we assume
$
(u,m)\neq(v,n)$ and $ m,n\ge1.$

Fix $\theta\in[0,1)$. Observe that
\begin{align*}
C:&=
\left\{
\rho^{d((u,m),(v,n))}
\bb{1}_{\{|(u,m)\wedge(v,n)|>0\}}
>\theta
\right\}\\[0.25cm]
&=
\left\{
d((u,m),(v,n))<K_\theta,\,
|(u,m)\wedge(v,n)|>0
\right\},
\end{align*}
where
$
K_\theta:=\log\theta/\log\rho,
$
with the convention that $\log0=-\infty$.

Thus, $C$ is the event that $(u,m)$ and $(v,n)$ belong to the same tree and their lowest common ancestor
$
L=(w,i)
$
satisfies $i>0$ together with
$$
d((u,m),L)+d((v,n),L)<K_\theta.
$$

Again, two cases arise.
\vspace{0.25cm}

\noindent
\textit{Case 1.} $(u,m)$ and $(v,n)$ belong to different branches below $L$.

In this case, the forest contains a non-root configuration shown in \Cref{fig:distance_case1}, where $i>0$ and $p+q<K_\theta$.

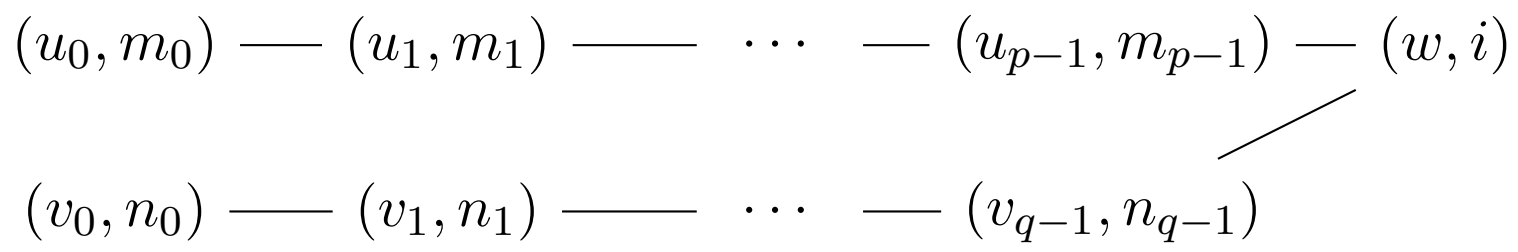
\begin{figure}[h]
\centering
\begin{tikzpicture}[scale=1.1]
\tikzset{
edge/.style={thick},
every node/.style={font=\small}
}

\node (u1) at (0,1) {$(u_0,m_0)$};
\node (u2) at (2,1) {$(u_1,m_1)$};
\node (up) at (6,1) {$(u_{p-1},m_{p-1})$};

\node (v1) at (0,0) {$(v_0,n_0)$};
\node (v2) at (2,0) {$(v_1,n_1)$};
\node (vq) at (6,0) {$(v_{q-1},n_{q-1})$};

\node (w1) at (8,1) {$(w,i)$};

\draw (u1)--(u2);
\draw (u2)--(3.5,1);
\node at (4,1) {$\cdots$};
\draw (4.5,1)--(up);

\draw (v1)--(v2);
\draw (v2)--(3.5,0);
\node at (4,0) {$\cdots$};
\draw (4.5,0)--(vq);

\draw (up)--(w1);
\draw (vq)--(w1);

\end{tikzpicture}
\caption{A non-root configuration corresponding to the event
$d((u,m),(v,n))<K_\theta$ when $(u,m)$ is not an ancestor of $(v,n)$.}
\label{fig:distance_case1}
\end{figure}

\vspace{0.25cm}

\noindent
\textit{Case 2.} $(u,m)$ lies on the ancestral line of $(v,n)$.

In this case, the forest contains the non-root configuration shown in \Cref{fig:distance_case2}, where $p<K_\theta$.

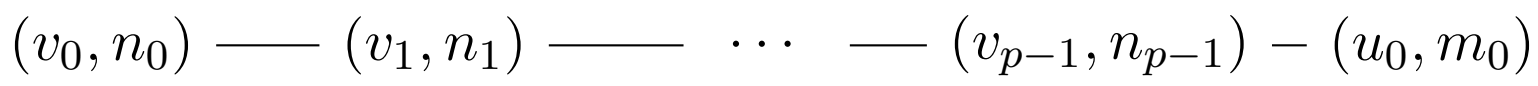
\begin{figure}[h]
\centering
\begin{tikzpicture}[scale=1.1]
\tikzset{
edge/.style={thick},
every node/.style={font=\small}
}

\node (u1) at (0,1) {$(v_0,n_0)$};
\node (u2) at (2,1) {$(v_1,n_1)$};
\node (up) at (6,1) {$(v_{p-1},n_{p-1})$};

\node (w1) at (8,1) {$(u_0,m_0)$};

\draw (u1)--(u2);
\draw (u2)--(3.5,1);
\node at (4,1) {$\cdots$};
\draw (4.5,1)--(up);

\draw (up)--(w1);

\end{tikzpicture}
\caption{A non-root configuration corresponding to the event
$d((u,m),(v,n))<K_\theta$ when $(u,m)$ is an ancestor of $(v,n)$.}
\label{fig:distance_case2}
\end{figure}

In both figures, the corresponding event is a non-root configuration. Since the initial vertices are fixed, these configurations are mutually disjoint. Hence $C$ is admissible. Applying \Cref{lem:admissible_comparison}, we obtain
$$
\Q_{\bF}(C)
\le
\Q_{\bF^1}(C^{(-1)}).
$$

Observe that $C^{(-1)}$ is precisely the event that $(u,m)$ and $(v,n)$ belong to the same tree and
$
d((u,m),(v,n))<K_\theta.
$
Equivalently,
\begin{align*}
C^{(-1)}
&=
\left\{
d((u,m),(v,n))<K_\theta,\,
\R{(u,m)}=\R{(v,n)}
\right\}\\
&=
\left\{
\rho^{d((u,m),(v,n))}
\bb{1}_{\{\R{(u,m)}=\R{(v,n)}\}}
>\theta
\right\}.
\end{align*}
This proves \eqref{eq:distance_compare_general} and concludes the proof.
\end{proof}

\newcommand{\iL}[2]{L_{#1}^{(#2)}}

Next we establish connection of these statistics with the same in a Hoppe tree. To proceed further we need two simple yet important concepts.  

\begin{definition}
An \emph{in-neighboring} (or \emph{backward}) \emph{random walk}
$
\{\iL{m}{v}:m\ge0\}
$
on $\G$, started from
$
\iL{0}{v}=v,
$
is the Markov chain with transition probabilities
$$
\Q\left(
\iL{m}{v}=u
\,\middle|\,
\iL{m-1}{v}=w
\right)
=
\begin{cases}
\displaystyle
\frac{1}{d_w^{\text{in}}}\,
\bb{1}_{\{u\in N_w\}},
&\text{if}~d_w^{\text{in}}>0,
\\[0.25cm]
\bb{1}_{\{u=w\}},
&\text{if}~d_w^{\text{in}}=0
\end{cases}
$$
for every $m\ge1$.
\end{definition}

\begin{definition}
Let $\bF$ be any rooted forest. The \emph{backward ancestral chain}
$
\{X_l:l\ge0\}
$ in $\bF$
starting from a vertex $x$ is defined recursively by $X_0=x$, and for each $l\ge1$,
$$
X_l=
\begin{cases}
\text{the parent of }X_{l-1},
& \text{if }X_{l-1}\text{ is not a root},
\\[0.2cm]
X_{l-1},
& \text{otherwise}.
\end{cases}
$$
\end{definition}

Clearly, the depth $|x|_{\bF}$ is completely determined by the backward ancestral chain $\{X_l\}$. Furthermore, if $\{Y_l\}$ is the backward ancestral chain starting from another vertex $y$, then both the depth of the lowest common ancestor, $|x\wedge y|_{\bF}$, and the graph distance $d_{\bF}(x,y)$ are completely determined by the pair of chains $\{X_l\}$ and $\{Y_l\}$.

\begin{lemma}\label{lem:comparison_with_hoppe_tree}
Suppose that $d_v^{\text{in}}>0$ for every $v\in V$, and let
$\mathcal T^t$ be a Hoppe tree with parameter $t>0$.
Then, for every $v\in V$ and $m,n\ge0$, the following relations hold:
\begin{align}
|(v,n)|_{\bF^t}
&\eqd
|n|_{\mathcal T^t},
\label{eq:depth_distribution}
\\
|(v,m)\wedge(v,n)|_{\bF^t}
&\st
|m\wedge n|_{\mathcal T^t},
\label{eq:lca_domination}
\\
d_{\mathcal T^t}(m,n)
&\st
d_{\bF^t}((v,m),(v,n)).
\label{eq:distance_domination}
\end{align}
\end{lemma}

\begin{proof}
First observe that \eqref{eq:depth_distribution} is trivial when $n=0$. Likewise, \eqref{eq:lca_domination} is trivial if either $m=0$ or $n=0$. If $m=n$, then \eqref{eq:lca_domination} reduces to \eqref{eq:depth_distribution}, while \eqref{eq:distance_domination} is immediate. Finally, when $m=0$ and $n>0$, \eqref{eq:distance_domination} is equivalent to \eqref{eq:depth_distribution}. Therefore, it suffices to prove \eqref{eq:depth_distribution}--\eqref{eq:distance_domination} for $m,n\ge1$ with $m\neq n$.

Let $\{X_l\}$ and $\{Y_l\}$ be the backward ancestral chains in $\bF^t$ starting from $(v,m)$ and $(v,n)$, respectively. Write
$$
X_l=(V_l^X,T_l^X),
\hspace{0.9cm}
Y_l=(V_l^Y,T_l^Y).
$$
By the construction of $\bF^t$ (see \eqref{eq:def_Ivn} and \eqref{eq:def_J}), the vertex parts $\{V_l^X:l\ge0\}$ and $\{V_l^Y:l\ge0\}$ are in-neighboring random walks started from $v$. Moreover,
$
\{V_l^X,V_l^Y:l\ge0\}
$
is independent of the time parts
$
\{T_l^X,T_l^Y:l\ge0\}.
$

We first prove \eqref{eq:depth_distribution}. By definition,
\begin{equation}\label{eq:sss1}
|(v,n)|_{\bF^t}
=
\inf\{l\ge0:T_l^Y=0\}.
\end{equation}
Since the chain $\{T_l^Y\}$ has the same distribution as the backward ancestral chain in the Hoppe tree $\mathcal T^t$ started from $n$, \eqref{eq:depth_distribution} follows immediately.

We next prove \eqref{eq:lca_domination} and \eqref{eq:distance_domination}. Define
$$
\tau
=
\min
\left\{
l\ge0:
\{T_r^X:r\le l\}\cap\{T_r^Y:r\ge0\}\neq\emptyset
\right\},
$$
that is, $\tau$ is the first time that the trajectory of $\{T^X_l\}$ intersects the trajectory of $\{T^Y\}$. Define a process $\{T_l:l\ge0\}$ by
$$
T_l=
\begin{cases}
T_l^X,&l\le\tau,\\
T_l^Y,&l>\tau.
\end{cases}
$$
Thus, $\{T_l\}$ follows the chain $\{T^X_l\}$ until it first intersects the trajectory of $\{T^Y_l\}$, after which it follows $\{T^Y_l\}$.

It is straightforward to verify that the pair $\{(T_l,T_l^Y):l\ge0\}$ has the same joint distribution as two backward ancestral chains in the Hoppe tree $\mathcal T^t$, started from $m$ and $n$, respectively. Let
$$
L(m,n):=T_\tau^X
$$
denote their first meeting point, and define
$$
\tau'
:=
\inf\{l\ge0:T_l^Y=L(m,n)\}.
$$
Notice that $L(m,n)$ has the same distribution as the lowest common ancestor $m\wedge n$ in $\mathcal T^t$. Consequently,
$$
\inf\{l\ge0:T_{\tau+l}=0\}
=
\inf\{l\ge0:T_{\tau'+l}^Y=0\}
\eqd
|m\wedge n|_{\mathcal T^t}.
$$

Suppose that
$
|(v,m)\wedge(v,n)|_{\bF^t}>0.
$
Then the chains $\{X_l\}$ and $\{Y_l\}$ meet before either reaches a root, and their meeting point is precisely the LCA $(v,m)\wedge(v,n)$. If
$
(v,m)\wedge(v,n)=Y_p,
$
then the chains $\{X_l\}$ and $\{Y_l\}$ cannot meet before the trajectories of $\{T^X_l\}$ and $\{T^Y_l\}$ intersect. Hence $p\le\tau'$, which implies
$$
|(v,m)\wedge(v,n)|_{\bF^t}
\le
\inf\{l\ge0:T_{p+l}^Y=0\}
\le
\inf\{l\ge0:T_{\tau'+l}^Y=0\}.
$$
Therefore,
\begin{equation}\label{eq:sss}
|(v,m)\wedge(v,n)|_{\bF^t}
\le
\inf\{l\ge0:T_{\tau'+l}^Y=0\},
\end{equation}
and the same inequality trivially holds when
$
|(v,m)\wedge(v,n)|_{\bF^t}=0.
$
Since the right-hand side has the same distribution as
$
|m\wedge n|_{\mathcal T^t},
$
this proves \eqref{eq:lca_domination}.

Finally, define
$$
\tilde d(m,n)
=
\inf\{l\ge0:T_l^X=L(m,n)\}
+
\inf\{l\ge0:T_l^Y=L(m,n)\}.
$$
By construction,
$
\tilde d(m,n)
\eqd
d_{\mathcal T^t}(m,n).
$
Moreover,
\begin{align*}
\tilde d(m,n)
&=
\inf\{l\ge0:T_l^X=0\}
+
\inf\{l\ge0:T_l^Y=0\}
\\
&\hspace{1cm}-
2\inf\{l\ge0:T_{\tau'+l}^Y=0\}\\[0.25cm]
&\le
|(v,m)|_{\bF^t}
+
|(v,n)|_{\bF^t}
-
2|(v,m)\wedge(v,n)|_{\bF^t}\\[0.25cm]
&=
d_{\bF^t}((v,m),(v,n)),
\end{align*}
where the inequality follows from \eqref{eq:sss1} and \eqref{eq:sss}. This proves \eqref{eq:distance_domination}.
\end{proof}

\begin{corollary}\label{cor:properties_depth_LCA}
Suppose $\bF$ is an interacting random recursive forest associated with
$\G$, of the first, second, or third kind, with positive root weights
$\{c_v:v\in V\}$ satisfying
$$
c:=\sup_{v\in V}c_v<\infty.
$$
Assume that $d_v^{\text{in}}>0$ for every $v\in V$. Then, for every $v\in V$ and integer $K\ge1$, the following hold as $n\to\infty$:
\begin{align}
\frac{|(v,n)|_{\bF}}{\log n}
&\xrightarrow{\Q}
1,
\label{eq:depth_asymptotic}
\\
\limsup_{n>m\to\infty}
\Q\Big(
|(v,m)\wedge(v,n)|_{\bF}>K
\Big)
&\le
\left(\frac12\right)^{K-1}.
\label{eq:lca_tail_bound}
\end{align}
\end{corollary}

\begin{proof}
By \Cref{lem:comparison_with_IRRF2,lem:comparison_with_hoppe_tree},
$$
|n|_{\mathcal T^c}
\st
|(v,n)|_{\bF}
\st
|n-1|_{\mathcal T^1}+1.
$$
Therefore, \eqref{eq:depth_asymptotic} follows immediately from \Cref{lem:Hoppe_depth}.

Next, combining \Cref{lem:comparison_with_IRRF2,lem:comparison_with_hoppe_tree}, we obtain
$$
\Q\Big(
|(v,m)\wedge(v,n)|_{\bF}>K
\Big)
\le
\Q\Big(
|(m-1)\wedge(n-1)|_{\mathcal T^1}>K-1
\Big).
$$
The right-hand side converges to
$$
\left(\frac12\right)^{K-1}
$$
as $n>m\to\infty$ by Theorem~7 of \cite{Lcaurt}. This proves
\eqref{eq:lca_tail_bound}.
\end{proof}

\begin{remark}
\textnormal{\eqref{eq:lca_tail_bound} in \Cref{cor:properties_depth_LCA} shows
that the lowest common ancestor of $(v,m)$ and $(v,n)$ remains at a
stochastically bounded distance from the root vertices in $\bF$. For Hoppe trees, this
stochastic boundedness of the lowest common ancestor (\Cref{lem:tightness_LCA_Hoppe_tree}) is the key ingredient in
the proof of Theorem~3.4 of \cite{bandyopadhyay2022new}, where convergence in
probability of the urn configurations is established. The same argument can
be adapted to obtain convergence in probability for interacting urn schemes
in greater generality (see \Cref{thm:IMVU_convergence}).
However, our goal in this section is to establish almost sure convergence,
for which substantially stronger estimates are required.}
\end{remark}

\begin{lemma}\label{lem:power_distance_comparison}
Suppose $\bF$ is an interacting random recursive forest associated with
$\G$, of the first, second, or third kind, with positive root weights
$\{c_v:v\in V\}$ satisfying
$$
c:=\sup_{v\in V}c_v<\infty.
$$
Assume that $d_v^{\text{in}}>0$ for every $v\in V$. Then, for every $v\in V$, $m,n\ge1$, and $\rho\in(0,1)$,
\begin{align}
\Ex{\rho^{d_{\bF}((v,m),(v,n))}
\bb1_{\{\R{(v,m)}=\R{(v,n)}\}}}
&\le
\Ex{\rho^{d_{\mathcal T^1}(m-1,n-1)}}
\notag\\
&+
(c+2)^2
\Ex{\rho^{|m-1|_{\mathcal T^1}}}
\Ex{\rho^{|n-1|_{\mathcal T^1}}},
\label{eq:power_sum_distance_dominance}
\end{align}
and 
\begin{align}
\Ex{\rho^{|(v,m)|_{\bF}+|(v,n)|_{\bF}}}
&\le
\Ex{\rho^{d_{\mathcal T^1}(m-1,n-1)}}
\notag\\
&+
2(c+2)^2
\Ex{\rho^{|m-1|_{\mathcal T^1}}}
\Ex{\rho^{|n-1|_{\mathcal T^1}}}.
\label{eq:power_sum_distance_dominance2}
\end{align}
\end{lemma}
\begin{proof}
Without loss of generality, assume that $m\neq n$. We first decompose the expectation on the left-hand side of \eqref{eq:power_sum_distance_dominance} as
\begin{align}
&\Ex{\rho^{d_{\bF}((v,m),(v,n))}
\bb1_{\{\R{(v,m)}=\R{(v,n)}\}}}\\
=&\Ex{\rho^{d_{\bF}((v,m),(v,n))}
\bb1_{\{|(v,m)\wedge(v,n)|_{\bF}>0\}}}
\notag\\
&+
\Ex{\rho^{d_{\bF}((v,m),(v,n))}
\bb1_{\{\R{(v,m)}=\R{(v,n)},\,|(v,m)\wedge(v,n)|_{\bF}=0\}}}.
\label{eq:www1}
\end{align}

For the first expectation on the right-hand side, \eqref{eq:distance_compare_general} of \Cref{lem:comparison_with_IRRF2} gives
\begin{align*}
&\Ex{\rho^{d_{\bF}((v,m),(v,n))}
\bb1_{\{|(v,m)\wedge(v,n)|_{\bF}>0\}}}\\
\le&
\Ex{\rho^{d_{\bF^1}((v,m-1),(v,n-1))}
\bb1_{\{\R{(v,m)}=\R{(v,n)}\}}}
\\
\le&
\Ex{\rho^{d_{\bF^1}((v,m-1),(v,n-1))}}.
\end{align*}
Now \eqref{eq:distance_domination} of \Cref{lem:comparison_with_hoppe_tree} implies that
$$
\Ex{\rho^{d_{\bF^1}((v,m-1),(v,n-1))}}
\le
\Ex{\rho^{d_{\T{}{1}}(m-1,n-1)}}.
$$
Hence,
\begin{equation}\label{eq:www}
\Ex{\rho^{d_{\bF}((v,m),(v,n))}
\bb1_{\{|(v,m)\wedge(v,n)|_{\bF}>0\}}}
\le
\Ex{\rho^{d_{\T{}{1}}(m-1,n-1)}}.
\end{equation}

Next, let
$$
G:=
\{\R{(v,m)}=\R{(v,n)},\,
|(v,m)\wedge(v,n)|_{\bF}=0\},
$$
which is the event that $(v,m)$ and $(v,n)$ belong to two different branches of the same tree of $\bF$. Consequently,
$$
\Ex{\rho^{d_{\bF}((v,m),(v,n))}
\bb1_G}
=
\Ex{\rho^{|(v,m)|_{\bF}+|(v,n)|_{\bF}}
\bb1_G}.
$$

Fix $i,j\ge1$. Let $u_1=v_1=v$, $m_1=m$, and $n_1=n$. Observe that the event
$$
\Omega_{i,j}
:=
\{|(v,m)|_{\bF}=i,\,
|(v,n)|_{\bF}=j\}
\cap G
$$
corresponds to configurations of the form:
\begin{figure}[H]
    \centering
    \begin{tikzpicture}[baseline=(current bounding box.center)]
        \node (u1) at (0,0) {$(u_1,m_1)$};
        \node (u2) at (2.2,0) {$(u_2,m_2)$};
        \node (udots) at (4.2,0) {$\cdots$};
        \node (ui) at (6.2,0) {$(u_i,m_i)$};
        \node (w) at (8.5,0) {$(w,0)$};

        \node (v1) at (0,-0.8) {$(v_1,n_1)$};
        \node (v2) at (2.2,-0.8) {$(v_2,n_2)$};
        \node (vdots) at (4.2,-0.8) {$\cdots$};
        \node (vj) at (6.2,-0.8) {$(v_j,n_j)$};

        \draw (u1) -- (u2);
        \draw (u2) -- (udots);
        \draw (udots) -- (ui);
        \draw (ui) -- (w);

        \draw (v1) -- (v2);
        \draw (v2) -- (vdots);
        \draw (vdots) -- (vj);
        \draw (vj) -- (w);
    \end{tikzpicture}
    \caption{}\label{fig:common_root_branches}
\end{figure}
where, for $l=2,\ldots,i$ and $s=2,\ldots,j$,
\begin{align}
u_l &\in N_{u_{l-1}},
\hspace{0.9cm}
v_s\in N_{v_{s-1}},
\label{eq:xxx1}
\\[0.25cm]
w &\in
\begin{cases}
N_{u_i}\cap N_{v_j},
    & \text{if } m_i,n_j>1,\\
\{u_i\},
    & \text{if } m_i=1,\ n_j>1,\\
\{v_j\},
    & \text{if } m_i>1,\ n_j=1,
\end{cases}
\label{eq:xyz1}
\\[-0.3ex]
&\hspace{2cm}\text{if $\bF$ is of first or second kind},
\nonumber
\\[0.5ex]
w &\in N_{v_i}\cap N_{v_j},
\hspace{0.9cm}
\text{if $\bF$ is of third kind},
\label{eq:xxx2}
\\[0.25cm]
m_{l-1} &> m_l,
\hspace{0.9cm}
n_{s-1}>n_s,
\label{eq:xxx3}
\\[0.25cm]
m_i,n_j &\geq 1,
\label{eq:xxx4}
\\[0.25cm]
\{(u_1,m_1),\ldots,(u_i,m_i)\}
&\cap
\{(v_1,n_1),\ldots,(v_j,n_j)\}
=\emptyset.
\label{eq:xxx5}
\end{align}
 Here, $(u_i,m_i)$ and $(v_j,n_j)$ are the vertices immediately preceding the common root on the ancestral paths from $(v,m)$ and $(v,n)$, respectively. Notice that $m_i$ and $n_j$ cannot both be equal to $1$ simultaneously when $\bF$ is of first or second kind.

Let $I$ denote the collection of all tuples
$$
\{(u_l,m_l),(v_s,n_s),w:l\leq i,\ s\leq j\}
$$
satisfying \eqref{eq:xxx1}--\eqref{eq:xxx5}, except for
\eqref{eq:xxx2} when $\bF$ is of the first or second kind, and except for
\eqref{eq:xyz1} when $\bF$ is of the third kind. For a fixed tuple
$$
\alpha=
\{(u_l,m_l),(v_s,n_s),w:l\le i,s\le j\}
\in I,
$$
let $A_\alpha$ denote the event that $\bF$ contains the configuration in \autoref{fig:common_root_branches}.

Assume first that the forest $\bF$ is of the first kind. Then, by definition,
\begin{align*}
\Q(A_\alpha)
=&
\left(
\prod_{l=1}^{i-1}
\frac{1}{c_{u_l}+(m_l-1)d_{u_l}^{\text{in}}}
\cdot
\frac{c_{u_i}}
{c_{u_i}+(m_i-1)d_{u_i}^{\text{in}}}
\right)\\
&\times
\left(
\prod_{s=1}^{j-1}
\frac{1}{c_{v_s}+(n_s-1)d_{v_s}^{\text{in}}}
\cdot
\frac{c_{v_j}}
{c_{v_j}+(n_j-1)d_{v_j}^{\text{in}}}
\right),
\end{align*}
where each product is interpreted as $1$ when $i=1$ or $j=1$, respectively.

Suppose first that $m_i,n_j>1$. Then clearly
\begin{equation}\label{bound1}
\Q(A_\alpha)
\le
c^2
\prod_{l=1}^i
\frac{1}{(m_l-1)d_{u_l}^{\text{in}}}
\prod_{s=1}^j
\frac{1}{(n_s-1)d_{v_s}^{\text{in}}}.
\end{equation}

Since $(u_1,m_1)=(v,m)$ and $(v_1,n_1)=(v,n)$ are fixed, the events
$\{A_\alpha:\alpha\in I\}$
are mutually disjoint. Therefore,
\begin{equation}\label{bound2}
\Q\left(\Omega_{i,j}\cap\{m_i,n_j>1\}\right)
\le
c^2
\sum_{\alpha\in I,\,
m_i,n_j>1}
\prod_{l=1}^i
\frac{1}{(m_l-1)d_{u_l}^{\text{in}}}
\prod_{s=1}^j
\frac{1}{(n_s-1)d_{v_s}^{\text{in}}}.
\end{equation}

If we omit the restrictions
\eqref{eq:xxx2} and \eqref{eq:xxx5}, the resulting sum only increases. Hence,
\begin{equation}\label{bound3}
\Q\left(\Omega_{i,j}\cap\{m_i,n_j>1\}\right)
\le
c^2
\sum_{\beta\in J,\,
m_i,n_j>1}
\prod_{l=1}^i
\frac{1}{(m_l-1)d_{u_l}^{\text{in}}}
\prod_{s=1}^j
\frac{1}{(n_s-1)d_{v_s}^{\text{in}}},
\end{equation}
where $J$ consists of all
$$
\beta=
\{(u_l,m_l),(v_s,n_s),w,x:l\le i,s\le j\}
$$
satisfying \eqref{eq:xxx1}, \eqref{eq:xxx3}, together with
$w\in N_{u_i}$ and $x\in N_{v_j}$.

Now summing over all choices of $u_l$ and $v_s$ satisfying
\eqref{eq:xxx1}, we obtain
$$
c^2
\sum_{m_1>m_2>\cdots>m_i>1}
\;
\sum_{n_1>n_2>\cdots>n_j>1}
\;
\prod_{l=1}^i
\frac{1}{m_l-1}
\prod_{s=1}^j
\frac{1}{n_s-1},
$$
which is precisely
$$
c^2
\Qr{|m-1|_{\T{}{1}}=i}
\Qr{|n-1|_{\T{}{1}}=j}.
$$
Hence,
$$
\Q\left(\Omega_{i,j}\cap\{m_i,n_j>1\}\right)
\le
c^2
\Qr{|m-1|_{\T{}{1}}=i}
\Qr{|n-1|_{\T{}{1}}=j}.
$$

Proceeding in the same way, we obtain
\begin{align*}
\Q\left(\Omega_{i,j}\cap\{m_i=1,n_j>1\}\right)
&\le c\,
\Qr{|m-1|_{\T{}{1}}=i-1}
\Qr{|n-1|_{\T{}{1}}=j},
\\[0.2cm]
\Q\left(\Omega_{i,j}\cap\{m_i>1,n_j=1\}\right)
&\le
c\,
\Qr{|m-1|_{\T{}{1}}=i}
\Qr{|n-1|_{\T{}{1}}=j-1}.
\end{align*}

Combining these estimates yields
\begin{align*}
\Q(\Omega_{i,j})
\le\;&
c^2
\Qr{|m-1|_{\T{}{1}}=i}
\Qr{|n-1|_{\T{}{1}}=j}
\\
&+
c\,
\Qr{|m-1|_{\T{}{1}}=i-1}
\Qr{|n-1|_{\T{}{1}}=j}
\\
&+
c\,
\Qr{|m-1|_{\T{}{1}}=i}
\Qr{|n-1|_{\T{}{1}}=j-1}.
\end{align*}

 The same estimate also holds when $\bF$ is of the second kind. However,
when $\bF$ is of the third kind, due to the additional case
$m_i=1,n_j=1$, we have one extra probability, which can be bounded as
$$
\Qr{\Omega_{i,j}\cap\{m_i=1,n_j=1\}}
\leq
\Qr{|m-1|_{\T{}{1}}=i-1}
\Qr{|n-1|_{\T{}{1}}=j-1}.
$$
Nevertheless, for all three kinds of $\bF$, we have
\begin{align*}
\Q(\Omega_{i,j})
\leq\;&
c^2
\Q\left(|m-1|_{\T{}{1}}=i\right)
\Q\left(|n-1|_{\T{}{1}}=j\right)
\\
&+
c\,
\Q\left(|m-1|_{\T{}{1}}=i-1\right)
\Q\left(|n-1|_{\T{}{1}}=j\right)
\\
&+
c\,
\Q\left(|m-1|_{\T{}{1}}=i\right)
\Q\left(|n-1|_{\T{}{1}}=j-1\right)
\\
&+
\Qr{|m-1|_{\T{}{1}}=i-1}
\Qr{|n-1|_{\T{}{1}}=j-1}.
\end{align*}

We now estimate
$$
\Ex{\rho^{|(v,m)|_{\bF}+|(v,n)|_{\bF}}
\bb1_G}
=
\sum_{i,j\ge1}
\rho^{i+j}
\Q(\Omega_{i,j}),
$$
which is bounded above by the sum of four terms. The first term is
\begin{align*}
&c^2
\sum_{i,j\ge1}
\rho^{i+j}
\Qr{|m-1|_{\T{}{1}}=i}
\Qr{|n-1|_{\T{}{1}}=j}\\
=&
c^2
\Ex{\rho^{|m-1|_{\T{}{1}}}}
\Ex{\rho^{|n-1|_{\T{}{1}}}}.
\end{align*}

The second term is
\begin{align*}
&
c
\sum_{i,j\ge1}
\rho^{i+j}
\Qr{|m-1|_{\T{}{1}}=i-1}
\Qr{|n-1|_{\T{}{1}}=j}
\\
&=
c
\left(
\rho\bb1_{\{m=1\}}
+
\Ex{\rho^{|m-1|_{\T{}{1}}}}
\right)
\Ex{\rho^{|n-1|_{\T{}{1}}}}
\\
&\le
2c\,
\Ex{\rho^{|m-1|_{\T{}{1}}}}
\Ex{\rho^{|n-1|_{\T{}{1}}}}.
\end{align*}

Similarly, the third term is bounded by
$$
2c\,
\Ex{\rho^{|m-1|_{\T{}{1}}}}
\Ex{\rho^{|n-1|_{\T{}{1}}}},
$$
while the fourth term is bounded by
$$
4\,
\Ex{\rho^{|m-1|_{\T{}{1}}}}
\Ex{\rho^{|n-1|_{\T{}{1}}}}.
$$

Combining these estimates, we obtain
\begin{equation}\label{eq:bound_disjointbranch_sameroot}
\Ex{\rho^{|(v,m)|_{\bF}+|(v,n)|_{\bF}}
\bb1_G}
\le
(c+2)^2
\Ex{\rho^{|m-1|_{\T{}{1}}}}
\Ex{\rho^{|n-1|_{\T{}{1}}}}.
\end{equation}
Together with \eqref{eq:www}, this proves
\eqref{eq:power_sum_distance_dominance}.

Finally, we now prove \eqref{eq:power_sum_distance_dominance2}. Let
$$
H:=\{\R{(v,m)}\neq\R{(v,n)}\},
$$
and define
$$
\Xi_{i,j}
:=
\{|(v,m)|_{\bF}=i,\,
|(v,n)|_{\bF}=j\}
\cap H.
$$
Thus, the event $\Xi_{i,j}$ occurs if and only if
$(v,m)$ and $(v,n)$ belong to different trees of $\bF$ with
$|(v,m)|_{\bF}=i$ and $|(v,n)|_{\bF}=j$. Equivalently, $\bF$
contains a configuration of the form:
\begin{figure}[H]
    \centering
    $$
    \begin{aligned}
    &(u_1,m_1)\longminus (u_2,m_2)\longminus \cdots
    \longminus (u_i,m_i)\longminus (w,0)\\[1ex]
    &(v_1,n_1)\longminus (v_2,n_2)\longminus \cdots
    \longminus (v_j,n_j)\longminus (x,0)
    \end{aligned}
    $$
    \caption{}\label{fig:two_branches_distinct_roots}
\end{figure}
where, for $l=2,\ldots,i$ and $s=2,\ldots,j$,
\begin{align}
u_l &\in N_{u_{l-1}},
\hspace{0.9cm}
v_s\in N_{v_{s-1}},
\label{eq:xxxx1}
\\[0.25cm]
(w,x)
&=(u_i,v_j),
\nonumber
\\[-0.1cm]
&\hspace{2cm}\text{if $\bF$ is of first or second kind},
\nonumber
\\[0.1cm]
(w,x)
&\in N_{u_i}\times N_{v_j},
\nonumber
\\[-0.1cm]
&\hspace{2cm}\text{if $\bF$ is of third kind},
\label{eq:xxxx2}
\\[0.25cm]
m_{l-1} &> m_l,
\hspace{0.9cm}
n_{s-1}>n_s,
\label{eq:xxxx3}
\\[0.25cm]
m_i,n_j &\geq 1,
\label{eq:xxxx4}
\\[0.25cm]
\{(u_1,m_1),\ldots,(u_i,m_i),(w,0)\}
&\cap
\{(v_1,n_1),\ldots,(v_j,n_j),(x,0)\}
=\emptyset.
\label{eq:xxxx5}
\end{align}

Proceeding exactly as in the case of the event $G$, by considering the four cases:
$$
\{m_i,n_j>1\},
\{m_i=1,n_j>1\},
\{m_i>1,n_j=1\},
\text{ and }
\{m_i=n_j=1\},
$$
we obtain
\begin{align*}
\Q(\Xi_{i,j})
\le\;&
c^2
\Qr{|m-1|_{\T{}{1}}=i}
\Qr{|n-1|_{\T{}{1}}=j}
\\
&+
c\,
\Qr{|m-1|_{\T{}{1}}=i-1}
\Qr{|n-1|_{\T{}{1}}=j}
\\
&+
c\,
\Qr{|m-1|_{\T{}{1}}=i}
\Qr{|n-1|_{\T{}{1}}=j-1}
\\
&+
\Qr{|m-1|_{\T{}{1}}=i-1}
\Qr{|n-1|_{\T{}{1}}=j-1}.
\end{align*}
Consequently,
$$
\Ex{\rho^{|(v,m)|_{\bF}+|(v,n)|_{\bF}}
\bb1_H}
=
\sum_{i,j\ge1}
\rho^{i+j}
\Q(\Xi_{i,j}),
$$
and the same calculation as before yields
\begin{equation}\label{eq:bound_diffroot}
\Ex{\rho^{|(v,m)|_{\bF}+|(v,n)|_{\bF}}
\bb1_H}
\le
(c+2)^2
\Ex{\rho^{|m-1|_{\T{}{1}}}}
\Ex{\rho^{|n-1|_{\T{}{1}}}}.
\end{equation}

Finally,
\begin{align*}
\Ex{\rho^{|(v,m)|_{\bF}+|(v,n)|_{\bF}}}
&=\Ex{\rho^{|(v,m)|_{\bF}+|(v,n)|_{\bF}}
\bb1_{\{|(v,m)\wedge(v,n)|_{\bF}>0\}}}\\
&+
\Ex{\rho^{|(v,m)|_{\bF}+|(v,n)|_{\bF}}
\bb1_G}+\Ex{\rho^{|(v,m)|_{\bF}+|(v,n)|_{\bF}}
\bb1_H}
\\[0.25cm]
&\le
\Ex{\rho^{d_{\bF}((v,m),(v,n))}
\bb1_{\{|(v,m)\wedge(v,n)|_{\bF}>0\}}}
\\
&+\Ex{\rho^{|(v,m)|_{\bF}+|(v,n)|_{\bF}}
\bb1_G}+
\Ex{\rho^{|(v,m)|_{\bF}+|(v,n)|_{\bF}}
\bb1_H}.
\end{align*}
Applying the bounds
\eqref{eq:www},
\eqref{eq:bound_disjointbranch_sameroot},
and
\eqref{eq:bound_diffroot}
establishes
\eqref{eq:power_sum_distance_dominance2}, thereby completing the proof.
\end{proof}

\subsection{Convergence of urn proportions}\label{sec:proof_mainresult}
The main goal of this part of the section is to prove the following theorem.

\begin{theorem}\label{thm:main_result}
Let
$
\{\C{n}{v}:v\in V,\;n\ge0\}
$
be an interacting urn scheme with reinforcement-based interaction, defined in
\Cref{subsec:mainmodel}, on a locally finite interaction graph
$\G=(V,\E)$. Suppose that assumptions \ref{A1} and \ref{A2} hold. In
addition, assume that:
\begin{enumerate}[label=(3.\arabic*),start=1]
    \item\label{A3}
    every vertex of $\G$ has at least one in-neighbor, that is,
    $$
    d_v^{\text{in}}>0,
    \hspace{0.9cm} \forall v\in V;
    $$

    \item\label{A4}
    the initial total number of balls is uniformly bounded, that is,
    $$
    \sup_{v\in V}T_0^{(v)}<\infty;
    $$

    \item\label{A5}
    the reinforcement matrix $R$ is irreducible and aperiodic.
\end{enumerate}
Then, for every $v\in V$,
$$
\frac{\C{n}{v}}{T_n^{(v)}}
\xrightarrow{a.s.}
\pi,
\hspace{0.9cm}\text{as }n\to\infty,
$$
where $\pi$ is the unique probability vector satisfying
$
\pi R=r\pi.
$
Moreover,
$$
\mathcal N_n^{(v)}
\xrightarrow{a.s.}
\pi,
\hspace{0.9cm}\text{as }n\to\infty.
$$
\end{theorem}

\begin{remark}
\textnormal{The main novelty of \Cref{thm:main_result} is that it remains valid for infinite graphs $\G$, thereby extending the finite-graph result of \cite{kaur2023interacting}. We believe that the aperiodicity assumption is not necessary. However, proving the result without this assumption would require stronger asymptotic estimates than those obtained in this article, and hence we assume aperiodicity throughout.}
\end{remark}

\begin{remark}
\textnormal{Another important observation is that the convergence in
\Cref{thm:main_result} does not depend on the particular topology of the
underlying graph $\G$. As long as $\G$ is locally finite and every vertex
has positive in-degree, the urn proportions at every vertex converge to
$\pi$. In particular, our result also applies to undirected graphs by
viewing each undirected edge as a pair of oppositely directed edges.}
\end{remark}

We follow the approach of \cite{bandyopadhyay2020strong} for establishing the almost sure convergence in \Cref{thm:main_result}. The estimates for several graph statistics related to interacting random recursive forests established in \Cref{subsec:properties_of_IRRF_statistics} enable us to adapt their argument to our setting.

Before proceeding, we require one final technical lemma, which extends Lemma 3.1 of \cite{bandyopadhyay2020strong} to the interacting urn setting.

\begin{lemma}\label{lem:covariance_bound}
Suppose $\bF$ is an interacting random recursive forest associated with
$\G$, of the first, second, or third kind, with positive root weights
$\{c_v:v\in V\}$ satisfying
$$
c:=\sup_{v\in V}c_v<\infty.
$$
Assume that $d_v^{\text{in}}>0$ for every $v\in V$. Let
$
\{W_{(v,n)}:v\in V,\ n\ge1\}
$
be the branching Markov process on $\bF$ defined in
\Cref{subsec:branchingmarkovprocess}. Assume further that $R$ is balanced, irreducible, and aperiodic. Then there exist constants $C>0$ and $\rho\in(0,1)$, depending only on $R$, such that, for every $v\in V$ and $m,n\ge1$,
\begin{align*}
&\max_{1\le i,j\le k}
\left|
Cov\!\left(
\bb1_{\{W_{(v,m)}=i\}},
\bb1_{\{W_{(v,n)}=j\}}
\,\middle|\,\Fo{}
\right)
\right|\\
&\hspace{6cm}\le
C
\rho^{d_{\Fo{}}((v,m),(v,n))}
\bb1_{\{\R{(v,m)}=\R{(v,n)}\}}.
\end{align*}
\end{lemma}

\begin{proof}
Fix $i,j\in\{1,\ldots,k\}$. For $m=n$, the inequality is trivial for any $C\ge1$. Hence, without loss of generality, assume that $m\neq n$.

Let $P$ and $Q$ denote the paths in the trees containing $(v,m)$ and $(v,n)$, respectively, starting from their root vertices
$\R{(v,m)}$ and $\R{(v,n)}$. Conditioned on $\Fo{}$, the processes
$\{W_a:a\in P\}$ and $\{W_b:b\in Q\}$ are Markov chains with transition matrix
$
S:=R/r,
$
where $r$ is the common row sum of $R$.

If $\R{(v,m)}\neq\R{(v,n)}$, then the two chains are independent, and hence the conditional covariance is zero. Therefore, assume that
$\R{(v,m)}=\R{(v,n)}$. In this case, the two chains evolve together until they reach their lowest common ancestor
$$
L:=(v,m)\wedge(v,n),
$$
after which they evolve independently (conditional on $W_L$). Consequently,
\begin{align}
&Cov\!\left(
\bb 1_{\{W_{(v,m)}=i\}},
\bb 1_{\{W_{(v,n)}=j\}}
\,\middle|\,\Fo{},W_L
\right)\notag\\
=&
\bb1_{\{\R{(v,m)}=\R{(v,n)}\}}
Cov\!\left(
S^{d_{\Fo{}}((v,m),L)}_{W_L,i},
S^{d_{\Fo{}}((v,n),L)}_{W_L,j}
\,\middle|\,\Fo{}
\right),\label{eq:covariance1}
\end{align}
where, by definition, 
$S^{d_{\Fo{}}((v,m),L)}_{W_L,i}$ and
$S^{d_{\Fo{}}((v,n),L)}_{W_L,j}$
are the conditional probabilities that, starting from the state $W_L$ at the lowest common ancestor, the two chains reach states $i$ and $j$ at $(v,m)$ and $(v,n)$, respectively.

To avoid notational complicacy we define $$d_{\bF}((v,n),L)=|(v,n)|_{\bF}-|(u,m)\wedge(v,n)|_{\bF}$$ and 
$$d_{\bF}((u,m),L)=|(u,m)|_{\bF}-|(u,m)\wedge(v,n)|_{\bF}
$$
even when $\R{(u,m)}\neq\R{(v,n)}$.

Now, observe that for any random variables $X$ and $Y$, and any constants $a,b$,
\begin{equation}\label{eq:xwxw1}
Cov(X,Y)
=
\Ex{X(Y-b)}
+
\Ex{X}\Ex{b-Y}.
\end{equation}
Choosing
$$
X=S^{d_{\Fo{}}((v,m),L)}_{W_L,i},
\hspace{0.9cm}
Y=S^{d_{\Fo{}}((v,n),L)}_{W_L,j},
$$
and
$$
a=\pi_i,
\hspace{0.9cm}
b=\pi_j,
$$
where
$\pi=(\pi_1,\ldots,\pi_k)$
is the stationary distribution of the Markov chain with transition matrix $S$, we obtain
$$
\left|
Cov\!\left(
S^{d_{\Fo{}}((v,m),L)}_{W_L,i},
S^{d_{\Fo{}}((v,n),L)}_{W_L,j}
\,\middle|\,\Fo{}
\right)
\right|
\le
2
\Ex{
\left|
S^{d_{\Fo{}}((v,n),L)}_{W_L,j}
-\pi_j
\right|
\,\middle|\,
\Fo{}
}.
$$

Since $S$ is irreducible and aperiodic, standard Markov chain theory (see, e.g., \cite{mixingtimes}) implies that there exist constants $C_0>0$ and $\rho^2\in(0,1)$, depending only on $S$, such that
\begin{equation}\label{eq:geometric_ergodicity}
\max_{1\le i\le k}
|S^n_{i,j}-\pi_j|
\le
C_0\rho^{2n},
\hspace{0.9cm}\forall n\ge1.
\end{equation}
Therefore,
$$
\left|
Cov\!\left(
S^{d_{\Fo{}}((v,m),L)}_{W_L,i},
S^{d_{\Fo{}}((v,n),L)}_{W_L,j}
\,\middle|\,\Fo{}
\right)
\right|
\le
2C_0
\rho^{2d_{\Fo{}}((v,n),L)}.
$$
By symmetry,
$$
\left|
Cov\!\left(
S^{d_{\Fo{}}((v,m),L)}_{W_L,i},
S^{d_{\Fo{}}((v,n),L)}_{W_L,j}
\,\middle|\,\Fo{}
\right)
\right|
\le
2C_0
\rho^{2d_{\Fo{}}((v,m),L)}.
$$
Combining the two estimates gives
\begin{align*}
\left|
Cov\!\left(
S^{d_{\Fo{}}((v,m),L)}_{W_L,i},
S^{d_{\Fo{}}((v,n),L)}_{W_L,j}
\,\middle|\,\Fo{}
\right)
\right|
&\le
2C_0
\rho^{2\max\{d_{\Fo{}}((v,m),L),\,d_{\Fo{}}((v,n),L)\}}\\
&\le
2C_0
\rho^{d_{\Fo{}}((v,m),(v,n))}.
\end{align*}
Substituting this bound into \eqref{eq:covariance1}, we obtain
$$
\left|
Cov\!\left(
\bb1_{\{W_{(v,m)}=i\}},
\bb1_{\{W_{(v,n)}=j\}}
\,\middle|\,\Fo{}
\right)
\right|
\le
2C_0
\rho^{d_{\Fo{}}((v,m),(v,n))}
\bb1_{\{\R{(v,m)}=\R{(v,n)}\}}.
$$
Finally, choosing $C\ge\max\{1,2C_0\}$ completes the proof.
\end{proof}

\noindent\textbf{Proof of \Cref{thm:main_result}.}
Let
$
S:=R/r.
$
Expanding the recursions \eqref{sum} and \eqref{mean}, we obtain
$$
\frac{\C{n}{v}-\C{0}{v}}{T_n^{(v)}}
=
\frac{\alpha_vrn}{T_n^{(v)}}
\frac1{d_v^{\text{in}}}
\sum_{u\in N_v}
\sum_{i=1}^k
S_i
\left(
\frac1n
\sum_{p=1}^n
\bb1_{\{\Z{p}{u}=i\}}
\right),
$$
where $\alpha_v=d_v^{\text{in}}$ for the sum-reinforcement scheme and
$\alpha_v=1$ for the average-reinforcement scheme.
Since
$$
T_n^{(v)}
=
T_0^{(v)}
+
\alpha_vrn
$$
and
$$
\sum_{i=1}^kS_i\pi_i=\pi,
$$
where $\pi=(\pi_1,\ldots,\pi_k)$, it therefore suffices to prove that
\begin{equation}\label{eq:colorcount_convergence}
\frac1n
\sum_{p=1}^n
\bb1_{\{\Z{p}{v}=i\}}
\xrightarrow{a.s.}
\pi_i,
\hspace{0.9cm}
\forall\,v\in V,\;
i=1,\ldots,k.
\end{equation}

Observe that
\begin{align}
\frac1n
\sum_{p=1}^n
\bb1_{\{\Z{p}{v}=i\}}
-\pi_i
=
&\frac1n
\sum_{p=1}^n
\left(
\bb1_{\{\Z{p}{v}=i\}}
-
\Qr{\Z{p}{v}=i}
\right)\notag\\
&+
\frac1n
\sum_{p=1}^n
\Qr{\Z{p}{v}=i}
-\pi_i.\label{eq:qqq}
\end{align}
We first prove that
\begin{equation}\label{eq:wqq}
\frac1n
\sum_{p=1}^n
\Qr{\Z{p}{v}=i}
-\pi_i
\xrightarrow{a.s.}
0.
\end{equation}
By \Cref{thm:representation_theorem}, it suffices to work with the branching Markov process
$
\{W_{(v,n)}:v\in V,\ n\ge1\}
$
on the interacting random recursive forest $\bF$ with interaction graph $\G$ and root weights
$c_v=T_0^{(v)}/r$, which is of the first or second kind according as we consider the sum- or average-reinforcement scheme.

Let $P$ denote the path in $\bF$ connecting the root $\R{(v,n)}$ to the
vertex $(v,n)$, and let $L_P$ denote the vertex of $P$ through which $P$ is
connected to its root $\R{(v,n)}$. Conditioned on $\bF$ and $W_{L_P}$, the event
$\{W_{(v,n)}=i\}$
is precisely the event that the Markov chain
$
\{W_a:a\in P\}
$
starting from $W_L$, reaches state $i$ after taking $|(v,n)|_{\bF}-1$ steps.
Therefore,
\begin{equation}
\label{eq:last}
\Qr{\W{v}{n}=i\,\middle|\,\bF,W_{L_P}}
=
S^{|(v,n)|_{\bF}-1}_{W_{L_P},i}.
\end{equation}
Since the chain is irreducible and aperiodic with stationary distribution $\pi$ and \eqref{eq:depth_asymptotic} of \Cref{cor:properties_depth_LCA} implies that
$
|(v,n)|_{\bF}
\xrightarrow{\Q}
\infty
$ as $n\rightarrow\infty$, by dominated convergence,
$$
\Qr{\W{v}{n}=i}\to\pi_i\hspace{0.9cm}\text{as}~n\to\infty.
$$
This proves \eqref{eq:wqq}. Hence, by \eqref{eq:qqq}, proving
\eqref{eq:colorcount_convergence} reduces to showing that
$$
\frac1n
\sum_{p=1}^n
\left(
\bb1_{\{\W{p}{v}=i\}}
-
\Qr{\W{p}{v}=i}
\right)
\xrightarrow{a.s.}
0.
$$

By Theorem~1 of \cite{SLLN_weakly_correlated}, it is enough to verify that
$$
\sum_{n=1}^{\infty}
\frac1{n^3}
Var
\left(
\sum_{p=1}^n
\bb1_{\{\W{v}{p}=i\}}
\right)
<
\infty,
$$
or equivalently,
\begin{equation}\label{eq:sum_cov_convergence_condition}
\sum_{N=1}^{\infty}
\frac1{N^3}
\sum_{1\le m\neq n\le N}
Cov
\left(
\bb1_{\{\W{v}{m}=i\}},
\bb1_{\{\W{v}{n}=i\}}
\right)
<
\infty.
\end{equation}

By the conditional covariance formula,
\begin{align}
Cov\!\left(
\bb1_{\{W_{(v,m)}=i\}},
\bb1_{\{W_{(v,n)}=i\}}
\right)
&=
\Ex{
Cov\!\left(
\bb1_{\{W_{(v,m)}=i\}},
\bb1_{\{W_{(v,n)}=i\}}
\,\middle|\,
\Fo{}
\right)
}
\notag
\\
&+
Cov\!\left(
\Q\!\left(
W_{(v,m)}=i
\,\middle|\,
\Fo{}
\right),
\Q\!\left(
W_{(v,n)}=i
\,\middle|\,
\Fo{}
\right)
\right).
\label{eq:xwxw2}
\end{align}
Let $J_{m,n}$ and $L_{m,n}$ denote the first and second terms on the right-hand side of \eqref{eq:xwxw2}, respectively.

Combining \Cref{lem:covariance_bound} with
\eqref{eq:power_sum_distance_dominance} of \Cref{lem:power_distance_comparison}, we obtain
$$
|J_{m,n}|
\le
\Ex{\rho^{d_{\T{}{1}}(m-1,n-1)}}
+
(c+2)^2
\Ex{\rho^{|m-1|_{\T{}{1}}}}
\Ex{\rho^{|n-1|_{\T{}{1}}}},
$$
where
$
c:=\sup_{v\in V}c_v
$
and $\rho\in(0,1)$ depends only on $R$ (as in \Cref{lem:covariance_bound}).
Now notice that \Cref{lem:power_distance_Hoppe_tree} yields
\begin{equation}\label{eq:asas}
\sum_{N=1}^{\infty}
\frac1{N^3}
\sum_{1\le m,n\le N}
\Ex{\rho^{d_{\T{}{1}}(m-1,n-1)}}
<
\infty,
\end{equation}
and
\begin{equation}\label{eq:bsbs}
\sum_{N=1}^{\infty}
\frac1{N^3}
\sum_{1\le m,n\le N}
\Ex{\rho^{|m-1|_{\T{}{1}}}}
\Ex{\rho^{|n-1|_{\T{}{1}}}}
\le
C
\sum_{N=1}^{\infty}
N^{2\rho-3}
<
\infty,
\end{equation}
for some constant $C>0$ depending only on $\rho$.
Consequently,
$$
\sum_{N=1}^{\infty}
\frac1{N^3}
\sum_{1\le m\neq n\le N}
|J_{m,n}|
<
\infty.
$$

Thus, to establish
\eqref{eq:sum_cov_convergence_condition}, it remains only to prove that
\begin{equation}\label{eq:cscs}
\sum_{N=1}^{\infty}
\frac1{N^3}
\sum_{1\le m\neq n\le N}
L_{m,n}
<
\infty.
\end{equation}

Observe that
\begin{align*}
\sum_{1\le m\neq n\le N}
L_{m,n}
&\le
Var
\left(
\sum_{n=1}^N
\Qr{\W{v}{n}=i\,\middle|\,\bF}
\right)
\\
&=
Var
\left(
\sum_{n=1}^N
\left(
\Qr{\W{v}{n}=i\,\middle|\,\bF}
-
\pi_i
\right)
\right)
\\
&\le
\Ex{
\left(
\sum_{n=1}^N
\left|
\Qr{\W{v}{n}=i\,\middle|\,\bF}
-
\pi_i
\right|
\right)^2
}.
\end{align*}

By the uniform geometric ergodicity estimate
\eqref{eq:geometric_ergodicity},
$$
\left|
\Qr{\W{v}{n}=i\,\middle|\,\bF}
-
\pi_i
\right|
\le
C_0
\rho^{2|(v,n)|_{\bF}}.
$$
Setting $\eta=\rho^2$, we obtain
\begin{align*}
\sum_{1\le m\neq n\le N}
L_{m,n}
&\le
C_0^2
\Ex{
\left(
\sum_{n=1}^N
\eta^{|(v,n)|_{\bF}}
\right)^2
}
\\
&=
C_0^2
\sum_{1\le m,n\le N}
\Ex{
\eta^{|(v,m)|_{\bF}+|(v,n)|_{\bF}}
}
\\
&\le
C_0^2
\sum_{1\le m,n\le N}
\Ex{\eta^{d_{\T{}{1}}(m-1,n-1)}}
\\
&\hspace{0.9cm}
+
2C_0^2(c+2)^2
\sum_{1\le m,n\le N}
\Ex{\eta^{|m-1|_{\T{}{1}}}}
\Ex{\eta^{|n-1|_{\T{}{1}}}},
\end{align*}
where the last inequality follows from
\eqref{eq:power_sum_distance_dominance2} of \Cref{lem:power_distance_comparison}. Since \eqref{eq:asas} and \eqref{eq:bsbs} hold for every $\rho\in(0,1)$, they imply that
\eqref{eq:cscs} holds. Consequently,
\eqref{eq:sum_cov_convergence_condition} follows, completing the proof. \qed

\subsection{Undirected graphs and self-organized criticality}
\label{subsec:SOC}

An undirected graph $\mathbb G:=(V,E)$ can be viewed as a directed graph in which
every edge is bidirectional, that is,
$(u,v)\in E$ implies $(v,u)\in E$ for all $u,v\in V$.
This includes, in particular, standard lattices such as
$\mathbb{Z}^d$, Bethe lattices, and Cayley graphs of finitely
generated groups, which arise naturally in the study of interacting
particle systems and statistical mechanics.
Recall that the convergence established in \Cref{thm:main_result}
does not depend on the particular topology of the underlying
directed graph, provided it is locally finite and every vertex
has positive in-degree. Consequently, the same result holds
for any locally finite undirected graph in which every vertex
has positive degree. In particular, the urn proportions at every
vertex converge almost surely to the same limiting vector $\pi$,
independently of the geometry of the graph.
We record this observation in the following corollary.

\begin{corollary}
\label{cor:undirected}
Let $\mathbb G=(V,E)$ be a locally finite undirected graph with no
isolated vertices. Suppose that the assumptions \ref{A1},\ref{A2},\ref{A4} and \ref{A5} are satisfied. Then, for every $v\in V$,
\[
\frac{\C{n}{v}}{T_n^{(v)}}\xrightarrow{{a.s.}}\pi,
\qquad \text{as } n\to\infty,
\]
where $\pi$ is the unique probability vector satisfying
$
\pi R=r\pi.
$
\end{corollary}

As mentioned in the introduction, one of our original motivations
for studying infinite interacting urns was to obtain a simple,
probabilistic, and exactly solvable model exhibiting behavior
reminiscent of self-organized criticality (SOC) (see
\cite{bak1987self,frigg2003self,dhar1993self,pruessner2012self}).
Undirected graphs, particularly lattices, provide a natural
setting for exploring this connection, since interactions in
many classical models of SOC take place between neighboring
sites of a lattice.

To illustrate this connection, think of each urn as a vertex of
$\G$ containing $k$ different types of sand, and let
\[
s_n^{(v)}:=\frac{\C{n}{v}}{T_n^{(v)}}
\]
denote the vector of proportions of the different types of sand
at vertex $v$ at time $n$. The reinforcement mechanism governed
by $R$ may then be interpreted as a rule describing how sand
spills over between neighboring vertices, with the spilling
occurring simultaneously throughout the system.

Under this interpretation, \Cref{cor:undirected} shows that,
irrespective of the initial distribution
$\{s_0^{(v)}:v\in V\}$, repeated spilling drives the proportions
of the different types of sand at every vertex towards the same
limiting vector $\pi$. In particular, on an undirected lattice,
local interactions between neighboring vertices lead to a
common asymptotic configuration throughout the infinite system,
without requiring any external parameter tuning.

Moreover, $\pi$ need not be an absorbing configuration. Even if
the proportions at every vertex coincide with $\pi$ at some
finite time, subsequent random spilling may move them away from
this configuration. Nevertheless, the dynamics continue to
drive the proportions towards $\pi$ in the long run. Thus,
the system continues to evolve while organizing itself around
the same limiting configuration.

An interesting feature of this phenomenon is that it does not
depend on the particular geometry of the interaction graph.
Indeed, the same convergence holds on arbitrary locally finite
undirected graphs with no isolated vertices, including regular
lattices, infinite trees, and more general irregular networks.
Thus, the emergence of a common limiting configuration is
driven by the local reinforcement mechanism rather than the
specific structure of the underlying graph.

Although this behavior is reminiscent of self-organized
criticality, the convergence result alone does not establish
critical behavior in the usual sense, such as scale-free
avalanches or power-law distributions. Nevertheless, it
provides a simple and exactly solvable framework in which
local stochastic interactions give rise to a robust form of
self-organization on infinite graphs.

\section{Some further applications of the representation approach}\label{sec:application}
\newcommand{\sS}{\mathcal{S}}
\newcommand{\mS}{\mathcal{M}(S)}
\newcommand{\pS}{\mathcal{P}(S)}
\newcommand{\hr}{\hat{R}}
This section is dedicated to show that the representation approach developed in the previous
section can be applied much more broadly. We shall see that the basic idea of
representing the colors of interacting urns through branching Markov processes
on a suitable interacting random recursive forest also works for other types of
interacting urn schemes and related reinforced stochastic processes. The exact
construction of the forest may change from one model to another, but the basic
idea remains essentially the same.

\subsection{Interacting measure-valued urns on infinite graphs}\label{subsec:IMVU_definition}
Measure-valued generalizations of urn models are well known (see
\cite{bandyopadhyay2017polya,Measure_valued_Polya_urn,bandyopadhyay2022new}).
We shall see that a similar generalization can be made for interacting urn
schemes. 

Let $\G=(V,\E)$ be the interaction graph introduced in \Cref{subsec:mainmodel}. Let $S$ be a Polish space, $\sS$ be its Borel sigma algebra, and let $\mS$ and $\pS$ denote respectively the spaces of finite measures and probability measures on $(S,\sS)$.

Suppose
$
R:S\times\sS\longrightarrow[0,\infty)
$
is a Markov type finite kernel, that is for every $s\in S$, the set-function
$R(s,\cdot)$ is a finite measure on $(S,\sS)$, and for every
$A\in\sS$, the function
$s\mapsto R(s,A)$ is $\sS/\mathcal{B}_{[0,\infty)]}$-measurable.

Let
$
\{\C{0}{v}:v\in V\}\subset\mS
$
be the initial urn configurations satisfying
$$
T_0^{(v)}:=\C{0}{v}(S)>0,
\hspace{0.9cm}\forall v\in V.
$$

We recursively construct random measures
$
\{\C{n}{v}:n\ge0,\;v\in V\}
$ together with $S$-valued random variables
$
\{\Z{n}{v}:n\ge1,\;v\in V\}
$
as follows.

For $n\ge1$, conditional on the sigma-algebra
$$
\sigma\!\left(\C{m}{u}:u\in V,\;0\le m\le n-1\right),
$$
the random variables
$
\{\Z{n}{v}:v\in V\}
$
are independent and satisfy
$$
\Q\left(\Z{n}{v}\in A\;\middle|\;
\C{m}{u}:u\in V,\;m\le n-1\right)
=
\frac{\C{n-1}{v}(A)}{\C{n-1}{v}(S)},\hspace{0.9cm}\forall A\in\sS.
$$
Equivalently, for every non-empty finite subset $V_0\subseteq V$ and measurable sets
$\{A_v:v\in V_0\}$,
$$
\Q\!\left(
\Z{n}{v}\in A_v,\;
v\in V_0
\;\middle|\;
\C{m}{u}:u\in V,\;m\le n-1
\right)
=
\prod_{v\in V_0}
\frac{\C{n-1}{v}(A_v)}
{\C{n-1}{v}(S)}.
$$

The urn configurations evolve according to

$$
\C{n}{v}
=
\begin{cases}
\displaystyle
\C{n-1}{v}
+
\sum_{u\in N_v}
R(\Z{n}{u},\cdot),
&
\text{if}~d_v^{\text{in}}>0,
\\[0.35cm]
\C{n-1}{v},
&
\text{if}~d_v^{\text{in}}=0,
\end{cases}
$$

which is the measure-valued analogue of the sum-reinforcement scheme \eqref{sum}.
Similarly, replacing the evolution rule by

$$
\C{n}{v}
=
\begin{cases}
\displaystyle
\C{n-1}{v}
+
\frac1{d_v^{\text{in}}}
\sum_{u\in N_v}
R(\Z{n}{u},\cdot),
&
\text{if}~d_v^{\text{in}}>0,
\\[0.35cm]
\C{n-1}{v},
&
\text{if}~d_v^{\text{in}}=0,
\end{cases}
$$

gives the measure-valued analogue of the average-reinforcement scheme \eqref{mean}.\\
Observe that the process remains well defined even when the interaction graph $\G$ is infinite, as long as it is locally finite.

The measure $\C{n}{v}$ is called the configuration of the urn at vertex
$v$ at time $n\ge0$, while the random variable $\Z{n}{v}$ denotes the color drawn from that urn at time $n-1\ge0$.
The kernel $R$ is called the reinforcement kernel. The process $\{\C{n}{v}:v\in V,n\ge0\}$ is called the \emph{interacting measure-valued urn scheme} on $\G$ with initial configuration $\{\C{0}{v}:v\in V\}$ and reinforcement kernel $R$.

Throughout this section, we assume that $R$ is balanced, that is, there exists
$r>0$ such that
$$
R(s,S)=r,
\hspace{0.9cm}
\forall s\in S.
$$
Consequently,
$$
T_n^{(v)}
:=
\C{n}{v}(S)
=
\begin{cases}
T_0^{(v)}+nrd_v^{\text{in}},
&
\text{for sum-reinforcement},
\\[0.2cm]
T_0^{(v)}+nr,
&
\text{for average-reinforcement}.
\end{cases}
$$

Observe that, when
$
S=\{1,\ldots,k\}
$
and $R(s,\cdot)$ is identified with the $s$-th row of a reinforcement matrix,
the above construction reduces to the finite-color interacting urn scheme
with reinforcement-based interaction considered in \Cref{subsec:mainmodel}.

\subsubsection{Branching Markov process for interacting measure-valued urns}\label{asderf}

Suppose
$$
\bF\sim\operatorname{IRRF}_1\left(\G,\{r^{-1}T_0^{(v)}\}\right)
\hspace{0.9cm}\text{or}\hspace{0.9cm}
\bF\sim\operatorname{IRRF}_2\left(\G,\{r^{-1}T_0^{(v)}\}\right).
$$
We now define a branching Markov process
$
\{\W{v}{n}:n\ge0,\;v\in V\}
$
on the forest $\bF$ using the kernel $R$ just as \Cref{subsec:branchingmarkovprocess}.

First, assign colors $W_{(v,1)}$ to the first generation vertices $\Fo{1}$ independently according to the distribution
$$
\Q\left(\W{v}{1}\in A\right)=\frac{\C{0}{v}(A)}{T_0^{(v)}},
\hspace{0.9cm}\forall A\in V.
$$

Now for $n\ge2$, suppose that the colors have been assigned to all vertices in $\Fo{n-1}$. Next construct the forest $\Fo{n}$. Conditional on the forest $\Fo{n}$ and the colors
$
\{\W{u}{m}:m\le n-1,\;u\in V\},
$
the random variables
$
\{\W{v}{n}:v\in V\}
$
are independent and, for every $A\in\sS$,
\begin{align*}
&\Q\!\left(
\W{v}{n}\in A
\,\middle|\,
\Fo{n},\,
\{\W{u}{m}:u\in V,\;m\le n-1\}
\right)\\
=&
\begin{cases}
r^{-1}R\left(W_{D_{(v,n)}},A\right),
& \text{if}~\D{v}{n}=(v,0),
\\[0.25cm]
\dfrac{\C{0}{v}(A)}{T_0^{(v)}},
& \text{if}~\D{v}{n}\neq(v,0).
\end{cases}
\end{align*}
Recall that $D_{(v,n)}$ is the parent of the vertex $(v,n)$ and lies in the forest $\Fo{n-1}$.

The following theorem extends \Cref{thm:representation_theorem} to interacting measure-valued
urns.

\begin{theorem}[Representation of interacting measure-valued urns]
\label{thm:representation_IMVU}
Suppose $\{\Z{n}{v}:v\in V,n\ge1\}$ is the color sequence associated with the interacting measure-valued urn scheme defined in \Cref{subsec:IMVU_definition}, Then the following distributional identity holds:
$$
\{\Z{n}{v}:v\in V,\;n\ge1\}
\eqd
\{\W{v}{n}:v\in V,\;n\ge1\}
$$
where $\{\W{v}{n}:v\in V,n\ge1\}$ is the branching Markov process on the interacting random recursive forest $\bF$ defined in \Cref{asderf}. The underlying forest $\bF$ is of first or second kind according to whether the process $\{\C{n}{v}:v\in V,n\ge0\}$ uses sum-reinforcement or average-reinforcement scheme.
\end{theorem}

\begin{proof}
The proof is entirely analogous to that of \Cref{thm:representation_theorem}. The only observation needed is that, for every $n\ge1$,
$$
\C{n}{v}=\C{0}{v},
\hspace{0.9cm}\text{whenever }d_v^{\text{in}}=0,
$$
whereas, for $d_v^{\text{in}}>0$,

\begin{equation}\label{eq:decomp_C}
\C{n}{v}
=
\C{0}{v}
+
\begin{cases}
\displaystyle
\sum_{j=1}^n
\sum_{u\in N_v}
R(\Z{j}{u},\cdot),
&
\text{for sum-reinforcement,}
\\[0.6cm]
\displaystyle
\frac{1}{d_v^{\text{in}}}
\sum_{j=1}^n
\sum_{u\in N_v}
R(\Z{j}{u},\cdot),
&
\text{for average-reinforcement. }
\end{cases}
\end{equation}
The remainder of the proof is similar to that of
\Cref{thm:representation_theorem}.
\end{proof}

Define $\hr=R/r$. Let
$\{\hr^n\}_{n\ge0}$ be the semigroup generated by the Markov kernel $\hr$.
That is for every bounded measurable function $f:S\to\mathbb{R}$,
$$
\hr f(s)=\int_S f(s')\,\hr(s,ds'),\hspace{0.9cm}\forall s\in S,
$$
and recursively,
$$
\hr^{n+1}f(s)
=
\int_S\hr^n f(s')\,\hr(s,ds'),\hspace{0.9cm}\forall s\in S.
$$

\begin{theorem}[Convergence of interacting measure-valued urns]
\label{thm:IMVU_convergence} Suppose $\{\C{n}{v}:v\in V,n\ge0\}$ is an interacting measure-valued urn scheme on a locally finite interaction graph $\G$ as defined in \Cref{subsec:IMVU_definition}. Assume that
\begin{enumerate}[label=(4.\arabic*),start=1]
    \item\label{B1}
    every vertex of $\G$ has at least one in-neighbor, that is,
    $$
    d_v^{\text{in}}>0,
    \hspace{0.9cm} \forall v\in V;
    $$

    \item\label{B2}
    the initial total amount of balls is uniformly bounded, that is,
    $$
    \sup_{v\in V}T_0^{(v)}<\infty;
    $$

    \item\label{B3} 
    there exists a probability $\pi\in\pS$ such that
\begin{equation}
\label{eq:convergence_condition_IMVU}
\hr^n(s,\cdot)\implies\pi,
\hspace{0.9cm}\forall s\in S.
\end{equation}
\end{enumerate}
Then, for every bounded continuous function $f:S\to\mathbb{R}$ and every $v\in V$,
$$
\frac{1}{T_n^{(v)}}
\int_S f\,d\C{n}{v}
\xrightarrow{\Q}
\int_S f\,d\pi,\hspace{0.9cm}\text{as}~n\to\infty.
$$
\end{theorem}

\begin{proof}
Using the decomposition \eqref{eq:decomp_C}, we obtain
$$
\frac{1}{T_n^{(v)}}
\left(
\int_S f\,d\C{n}{v}
-
\int_S f\,d\C{0}{v}
\right)
=
\frac{\alpha_vrn}{T_n^{(v)}}
\frac{1}{d_v^{\text{in}}}
\sum_{u\in N_v}
\left(
\frac1n
\sum_{p=1}^n
\hr f(\Z{p}{u})
\right),
$$
where $T_n^{(v)}=T_0^{(v)}+\alpha_vrn$, and
$\alpha_v$ is as defined in the proof of
\Cref{thm:main_result}.
Hence it is enough to prove that
$$
\frac1n
\sum_{p=1}^n
\hr f(\Z{p}{v})
\xrightarrow{\Q}
\pi f,
\hspace{0.9cm}
v\in V,
$$
where
$$
\pi f:=\int_Sf\,d\pi.
$$

To this end, observe that
\begin{equation}
\label{eq:variance_IMVU}
Var\left(
\frac1n
\sum_{p=1}^n
\hr f(\Z{p}{v})
\right)
=
o(1)
+
\frac{2}{n^2}
\sum_{1\le p<q\le n}
Cov\!\left(
\hr f(\Z{p}{v}),
\hr f(\Z{q}{v})
\right),
\end{equation}
where $o(1)\to0$ as $n\to\infty$. We shall prove that the variance in
\eqref{eq:variance_IMVU}
converges to zero as $n\to\infty$.

Since
$
|\hr f|\le\|f\|_\infty,
$
(where $\|\cdot\|_\infty$ denotes the supremum norm) every covariance term in \eqref{eq:variance_IMVU} is uniformly bounded.
Therefore, by dominated convergence, it is enough to show that
\begin{equation}
\label{eq:cov_goingto_0}
\lim_{q>p\to\infty}
Cov\!\left(
\hr f(\Z{p}{v}),
\hr f(\Z{q}{v})
\right)
=
0.
\end{equation}

By the representation theorem
(\Cref{thm:representation_IMVU}),
it suffices to work with the branching Markov process $\W{v}{n}$ instead of $\Z{n}{v}$.

Let $P$ and $Q$ denote the branches of $\bF$ containing the vertices
$(v,p)$ and $(v,q)$, respectively. Let $L_P$ and $L_Q$ denote the
vertices of $\bF$ adjacent to the corresponding root vertices on the
branches $P$ and $Q$, respectively. Conditioned on the forest $\bF$,
the processes $\{W_a:a\in P\}$ and $\{W_b:b\in Q\}$ are Markov chains
with transition kernel $R$, started from the colors $W_{L_P}$ and
$W_{L_Q}$, respectively.

Suppose first that $P$ and $Q$ belong to two different trees of $\bF$. Then conditioned of $\bF$, $W_{L_P}$ and $W_{L_Q}$, the two chains evolve independently by construction and therefore
$$
\begin{aligned}
&
\Ex{
\hr f(\W{v}{p})
\hr f(\W{v}{q})
\,\middle|\,
\bF,W_{L_P},W_{L_Q}
}
\bb1_{\{\R{(v,p)}\neq\R{(v,q)}\}}
\\
&\hspace{0.9cm}=
\hr ^{|(v,p)|_{\bF}-1}\hr
f(W_{L_P})
\,
\hr^{|(v,q)|_{\bF}-1}\hr
f(W_{L_Q})
\,
\bb1_{\{\R{(v,p)}\neq\R{(v,q)}\}}
\\
&\hspace{0.9cm}=
\hr ^{|(v,p)|_{\bF}}
f(W_{L_P})
\,
\hr^{|(v,q)|_{\bF}}
f(W_{L_Q})
\,
\bb1_{\{\R{(v,p)}\neq\R{(v,q)}\}}
\end{aligned}
$$
Now assumption \ref{B3} and the fact that
$
|(v,p)|_{\bF}
\xrightarrow{\Q}
\infty
$ as $p\to\infty$
(see \Cref{cor:properties_depth_LCA}), give
$$
\hr^{|(v,p)|_{\bF}}
f(W_{L_P})
\xrightarrow{\Q}
\pi f,
\hspace{0.9cm}
\text{as}~p\to\infty.
$$
Hence, by dominated convergence,
\begin{equation}
\label{eq:unequalroots_estimate}
\lim_{q>p\to\infty}
\Ex{
\left|
\Ex{
\hr f(\W{v}{p})
\hr f(\W{v}{q})
\,\middle|\,
\bF
}
-
(\pi f)^2
\right|
\bb1_{A_{-1}}
}
=
0,
\end{equation}
where
$$
A_{-1}
=
\{
\R{(v,p)}
\neq
\R{(v,q)}
\}.
$$
Define
$$
A_j
=
\left\{
\R{(v,p)}=\R{(v,q)},
~
|(v,p)\wedge(v,q)|_{\bF}=j
\right\},
\hspace{0.9cm}\text{for}~j\ge0.
$$
Conditioned on the forest $\bF$, if $P$ and $Q$ belong to the same tree, then the two chains coincide along the common segment $P\cap Q$ and evolve separately after reaching their lowest common ancestor
$
L=(v,p)\wedge(v,q).
$
Consequently, conditioned on $\bF$ and $W_L$, if $A_j$ occurs, then the descendant chains
$$
\{W_a:a\in (P\setminus Q)\cup{L}\}
\quad\text{and}\quad
\{W_b:b\in (Q\setminus P)\cup\{L\}\}
$$
evolve independently. Therefore,
\begin{align*}
&\Ex{
\hr f(\W{v}{p})
\hr f(\W{v}{q})
\,\middle|\,
\bF, W_L
}
\bb1_{A_j}\\
&\hspace{4cm}=
\hr^{|(v,p)|_{\bF}-j+1}
f(W_L)
\hr^{|(v,q)|_{\bF}-j+1}
f(W_L)
\bb1_{A_j}.
\end{align*}
Since $j$ is fixed and
$
|(v,p)|_{\bF},
~
|(v,q)|_{\bF}
\xrightarrow{\Q}\infty
$ (as $p,q\to\infty$),
assumption \ref{B3} again yields
$$
\hr^{|(v,p)|_{\bF}-j+1}
f(W_L)\hr^{|(v,q)|_{\bF}-j+1}
f(W_L)
\xrightarrow{\Q}
(\pi f)^2,\hspace{0.9cm}\text{as}~p,q\to\infty.
$$
Hence, by dominated convergence,
\begin{equation}
\label{eq:equalroots_estimate}
\lim_{q>p\to\infty}
\Ex{
\left|
\Ex{
\hr f(\W{v}{p})
\hr f(\W{v}{q})
\,\middle|\,
\bF
}
-
(\pi f)^2
\right|
\bb1_{A_j}
}
=
0.
\end{equation}

For convenience, write
$$
M_{p,q}
=
\hr f(\W{v}{p})
\hr f(\W{v}{q}).
$$
Then, for every integer $K\ge1$,
\begin{align*}
\left|
\Ex{M_{p,q}}
-
\sum_{j=-1}^{K}
\Ex{
\Ex{M_{p,q}\,\middle|\,\bF}
\bb1_{A_j}
}
\right|
&\le
\sum_{j>K}
\Ex{
\Ex{|M_{p,q}|\,\middle|\,\bF}
\bb1_{A_j}
}
\\
&\le
\|f\|_\infty^2
\Qr{
|(v,p)\wedge(v,q)|_{\bF}>K
}.
\end{align*}
By \Cref{cor:properties_depth_LCA},
$$
\limsup_{q>p\to\infty}
\left|
\Ex{M_{p,q}}
-
\sum_{j=-1}^{K}
\Ex{
\Ex{M_{p,q}\,\middle|\,\bF}
\bb1_{A_j}
}
\right|
\le
\frac{2\|f\|_\infty^2}{2^K}.
$$
Combining this with
\eqref{eq:unequalroots_estimate}
and
\eqref{eq:equalroots_estimate},
we obtain
$$
\limsup_{q>p\to\infty}
\left|
\Ex{M_{p,q}}
-
(\pi f)^2
\Qr{
|(v,p)\wedge(v,q)|_{\bF}\le K
}
\right|
\le
\frac{2\|f\|_\infty^2}{2^K}.
$$
Finally, letting $K\to\infty$ gives
\begin{equation}
\label{eq:sdsdfg}
\lim_{q>p\to\infty}
\Ex{
\hr f(\W{v}{p})
\hr f(\W{v}{q})
}
=
(\pi f)^2.
\end{equation}
We now turn to the individual expectations. Observe that
$$
\Ex{\hr f(\W{v}{p})\,\middle|\,\bF,W_{L_P}}
=
R^{|(v,p)|_{\bF}}f(W_{L_P}).
$$
which converges to $\pi f$ due to assumption \ref{B3} and hence by dominated convergence,
\begin{equation}\label{eq:individual_expectation_IMVU}
\lim_{p\to\infty}
\Ex{\hr f(\W{v}{p})}
=
\pi f.
\end{equation}
Combining this with \eqref{eq:sdsdfg} we obtain
$$
\lim_{q>p\to\infty}
Cov\!\left(
\hr f(\W{v}{p}),
\hr f(\W{v}{q})
\right)
=
0,
$$
which proves \eqref{eq:cov_goingto_0}. Consequently, the variance in
\eqref{eq:variance_IMVU} converges to zero as $n\to\infty$. Now the only thing remains to prove is that
$$
\lim_{n\to\infty}\frac1n
\sum_{p=1}^n
\Ex{\hr f(\W{v}{p})}
=
\pi f
$$
which follows directly from \eqref{eq:individual_expectation_IMVU}. This completes the proof.
\end{proof}

\begin{remark}
\textnormal{\Cref{thm:IMVU_convergence} extends the convergence result for finite-color
interacting urn schemes with reinforcement-based interaction in
\Cref{sec:proof_mainresult} to Polish space-valued colors. However, the
convergence is in probability rather than almost surely. In the finite-color
case, irreducibility and aperiodicity of $R$ give the uniform geometric
ergodicity property \eqref{eq:geometric_ergodicity} for the corresponding
Markov chain. This gives us stronger covariance estimates (\Cref{lem:covariance_bound}), which are needed for almost sure
convergence. The same estimates can also be obtained in the measure-valued
setting if the Markov chain with kernel $\hr$ satisfies a uniform geometric
ergodicity condition such as \eqref{eq:geometric_ergodicity}. This is also the
condition used in \cite{bandyopadhyay2020strong} to prove almost sure
convergence of urn proportions.}
\end{remark}

\subsection{Interacting urn schemes with drawing-based interaction}
\label{subsec:drawing_interaction}
In this section, we study another type of interacting urn scheme
that is well known in the literature.

Let $\G=(V,\E)$ be the interaction graph introduced in
\Cref{subsec:mainmodel}. Throughout this section, we assume that
$d_v^{\text{in}}>0$ for every $v\in V$. Suppose each edge $\e\in\E$ is assigned a non-negative weight
$w_{\e}$ satisfying
$$
\sum_{u\in N_v}w_{(u,v)}=1,
\hspace{0.9cm}
\forall v\in V.
$$
We refer to $\{w_{\e}:\e\in E\}$ as the collection of edge weights.

The urn configurations $\C{n}{v}$ are $k$-dimensional row vectors, with
$k\geq 2$, and the drawn colors $\Z{n}{v}$ take values in
$\{1,\ldots,k\}$, as in \Cref{subsec:mainmodel}. The initial configuration $\{\C{0}{v}:v\in V\}$ and the reinforcement matrix
$R$ is assumed to satisfy conditions
\ref{A1} and \ref{A2}. The dynamics evolves recursively as follows.

For $n\ge1$, conditioned on the history
$
\{\C{m}{u}:u\in V,\;m\le n-1\},
$
the random variables
$
\{\Z{n}{v}:v\in V\}
$
are independent. More precisely, for every finite non-empty subset
$V_0\subseteq V$ and every choice of colors
$i_v\in\{1,\ldots,k\}$ for all $v\in V_0$, the conditional distribution is given by 
\begin{align}
&\Q\!\left(
\Z{n}{v}=i_v,\;
v\in V_0
\,\middle|\,
\C{m}{u}:u\in V,\;m\le n-1
\right)\notag\\
=&
\prod_{v\in V_0}
\sum_{u\in N_v}
w_{(u,v)}
\frac{\C{n-1}{u}(i_v)}
{T_{n-1}^{(u)}}.\label{eq:pok}
\end{align}
After the colors are drawn, the urns are updated according to
$$
\C{n}{v}
=
\C{n-1}{v}
+
R_{\Z{n}{v}},
\hspace{0.9cm}
\forall v\in V.
$$
Since the reinforcement matrix is balanced, the total number of balls satisfies
$$
T_n^{(v)}
=
T_0^{(v)}
+
nr,
\hspace{0.9cm}
\forall v\in V.
$$
The resulting process
$
\{\C{n}{v}:v\in V,\;n\ge0\}
$
is called an \emph{interacting urn scheme with drawing-based interaction}
on the interaction graph $\G$, with initial configuration
$\{\C{0}{v}:v\in V\}$, reinforcement matrix $R$, and edge weights
$\{w_{\e}:\e\in E\}$. The process
$
\{\Z{n}{v}:v\in V,\;n\ge1\}
$
is called the \emph{color sequence} of the interacting urn scheme.
\begin{remark}
\textnormal{When the interaction graph $\G$ is finite and complete, this model was introduced and studied in~\cite{aletti2017interacting}. Assuming that the reinforcement matrix is irreducible, the authors established the first- and second-order asymptotic behavior of the joint urn proportions using stochastic approximation.}
\end{remark}

\subsubsection{Interacting random recursive forests with edge weights}

We now introduce a weighted version of the interacting random recursive forest of third kind defined in
\Cref{subsec:properties_of_IRRF_statistics}.

Let
$
\{c_v:v\in V\}\subset(0,\infty)
$
be a collection of root weights and let
$
\{w_{\e}:\e\in E\}
$
be the edge weights introduced in
\Cref{subsec:drawing_interaction}. Define a family of independent random variables
$$
\{I_{(v,n)},J_{(v,n)}:v\in V,\;n\ge0\}
$$
as follows.

For every $v\in V$ and $n\ge0$,
$$
\Q(I_{(v,n)}=u)
=
w_{(u,v)}
\bb1_{\{u\in N_v\}}.
$$
Further let
$
J_{(v,0)}=J_{(v,1)}=0,
$
and for $n\ge2$,
$$
\Q(J_{(v,n)}=j)
=
\begin{cases}
\displaystyle
\frac{c_v}{c_v+n-1},
&
\text{if}~j=0,
\\[0.35cm]
\displaystyle
\frac{1}{c_v+n-1},
&
\text{if}~1\le j\le n-1.
\end{cases}
$$
The parent variables are defined by
$$
D_{(v,n)}
=
\left(I_{(v,n)},J_{(v,n)}\right).
$$

The forest is constructed recursively as follows. Let
$$
\Fo{0}
=
\{(v,0):v\in V\}
$$
be the root vertices, and for $n\ge1$, having constructed $\Fo{n-1}$, obtain $\Fo{n}$ by attaching each vertex $(v,n)$ to its parent $D_{(v,n)}$. The resulting forest
$$
\bF
=
\bigcup_{n\ge0}\Fo{n}
$$
is called the \emph{interacting random recursive forest of third kind} associated with the interaction graph $\G$, root weights $\{c_v:v\in V\}$, and edge weights $\{w_{\e}:\e\in E\}$. In short we write
$$
\bF
\sim
\operatorname{IRRF}_3^{\,w}
\left(
\G,
\{c_v:v\in V\},
\{w_{\e}:\e\in E\}
\right).
$$

Notice that if $
w_{(u,v)}
=
1/d_v^{\text{in}}$
 for all $u\in N_v$ and $v\in V,
$
then
$$
\operatorname{IRRF}_3^{\,w}
\left(
\G,
\{c_v:v\in V\},
\{w_{\e}:\e\in E\}
\right)
=
\operatorname{IRRF}_3\left(\G,\{c_v:v\in V\}\right).
$$

\subsubsection{Branching Markov process for interacting urn scheme with drawing-based interaction}
\label{subsubsec:asasd}

Let
$$
\bF\sim
\operatorname{IRRF}_3^{\,w}
\left(
\G,
\{r^{-1}T_0^{(v)}:v\in V\},
\{w_{\e}:\e\in E\}
\right).
$$
The branching Markov process
$
\{\W{v}{n}:v\in V,\;n\ge1\}
$
on the forest $\bF$ is defined just like previous sections.

First, assign colors to the first generation vertices of $\bF$ independently by letting
$$
\Qr{\W{v}{1}=i}
=
\frac{\C{0}{v}(i)}{T_0^{(v)}},
\hspace{0.9cm}\forall i=1,\ldots,k,~
v\in V.
$$
Now suppose that the colors have been assigned to all vertices of $\Fo{n-1}$. Having constructed the forest $\Fo{n}$, assign colors to the vertices of the
$n$-th generation as follows.

Conditioned on the forest $\Fo{n}$ and the
previously assigned colors
$
\{\W{v}{m}:v\in V,\;m\le n-1\},
$
the random variables
$
\{\W{v}{n}:v\in V\}
$
are independent and satisfy
\begin{align*}
&\Q\!\left(
\W{v}{n}=i
\,\middle|\,
\Fo{n},
\{\W{u}{m}:u\in V,\;m\le n-1\}
\right)\\[0.25cm]
=&
\begin{cases}
\displaystyle
r^{-1}R_{W_{D_{(v,n)}},\,i},
&
\text{if }D_{(v,n)}~\text{is not a root},
\\[0.4cm]
\displaystyle
\frac{\C{0}{u}(i)}{T_0^{(u)}},
&
\text{if }D_{(v,n)}=(u,0)~\text{for some}~u\in V.
\end{cases}
\end{align*}
For all $i=1,\ldots,k$. 

A rigorous construction can also be carried out in
the same manner as in
\Cref{subsec:branchingmarkovprocess}.
\begin{theorem}[Representation of interacting urn scheme with drawing-based interaction]\label{thm:reprentation_3}
  Suppose $\{\Z{n}{v}:v\in V,n\ge1\}$ is a color sequence associated with the interacting urn scheme defined in \Cref{subsec:drawing_interaction}. Then 
  $$
  \{\Z{n}{v}:v\in V,n\ge1\}\eqd\{\W{v}{n}:v\in V,n\ge1\}
  $$
  where $\{\W{v}{n}:v\in V,n\ge1\}$ is the branching Markov process defined in \Cref{subsubsec:asasd}.
\end{theorem}

\begin{proof}
\newcommand{\sg}[1]{\mathcal{G}_{#1}}

First observe that
$$
\C{n}{v}
=
\C{0}{v}
+
\sum_{j=1}^n
R_{\Z{j}{v}},
\hspace{0.9cm}
\forall n\ge1.
$$
We proceed exactly as in the proof of
\Cref{thm:representation_theorem}. Define
$
H_{(v,0)}:=\C{0}{v},
$
and, for $n\ge1$,
\begin{equation}
\label{aspok}
H_{(v,n)}
:=
H_{(v,0)}
+
\sum_{j=1}^n
R_{\W{v}{j}}.
\end{equation}
Let
$
\sg{n}
:=
\sigma\!\left(
\{\W{v}{m}:v\in V,\;m\le n\}
\right).
$
By construction, it is immediate that
$$
\{\Z{1}{v},\C{0}{v}:v\in V\}
\eqd
\{\W{v}{1},H_{(v,0)}:v\in V\},
$$
and therefore from \eqref{aspok}
$$
\{\C{1}{v}:v\in V\}
\eqd
\{H_{(v,1)}:v\in V\}.
$$

Now assume that, for some $n\ge2$,
$$
\{\Z{m}{v},\C{m-1}{v}:v\in V,\;m\le n-1\}
\eqd
\{\W{v}{m},H_{(v,m-1)}:v\in V,\;m\le n-1\}.
$$
In particular, from \eqref{aspok}
$$
\{\C{n-1}{v}:v\in V\}
\eqd
\{H_{(v,n-1)}:v\in V\}.
$$

Now since conditioned on $\Fo{n}$ and $\sg{n-1}$, the random variables
$\{\W{v}{n}:v\in V\}$
are independent, we have, for every
$i\in\{1,\ldots,k\}$,
$$
\Q\!\left(
\W{v}{n}=i
\,\middle|\,
\Fo{n},\sg{n-1}
\right)
=
\frac{H_{(v,0)}(i)}{T_0^{(v)}}
\bb1_{\{\D{v}{n}\in\Fo{0}\}}
+
r^{-1}
R_{W_{\D{v}{n}},i}
\bb1_{\{\D{v}{n}\notin\Fo{0}\}}.
$$
Taking expectation with respect to the forest $\Fo{n}$, we obtain
\begin{align*}
\Q\!\left(
\W{v}{n}=i
\,\middle|\,
\sg{n-1}
\right)
&=
\sum_{u\in N_v}
w_{(u,v)}
\frac{H_{(u,0)}(i)}{T_0^{(u)}}
\frac{r^{-1}T_0^{(u)}}{r^{-1}T_0^{(u)}+n-1}\\
&\hspace{0.9cm}\hspace{0.9cm}+
\sum_{u\in N_v}
w_{(u,v)}
\sum_{j=1}^{n-1}
\frac{r^{-1}R_{\W{u}{j},i}}
{r^{-1}T_0^{(u)}+n-1}
\\[0.2cm]
&=
\sum_{u\in N_v}
w_{(u,v)}
\frac{H_{(u,n-1)}(i)}
{T_{n-1}^{(u)}}.
\end{align*}
Moreover by construction, conditioned on $\sg{n-1}$, the random variables
$\{\W{v}{n}:v\in V\}$
remain independent. Consequently, their conditional distribution agrees with that of the colors
$\{\Z{n}{v}:v\in V\}$
given by \eqref{eq:pok}. Therefore,
$$
\{\Z{m}{v},\C{m-1}{v}:v\in V,\;m\le n\}
\eqd
\{\W{v}{m},H_{(v,m-1)}:v\in V,\;m\le n\},
$$
which completes the induction. Hence,
$$
\{\Z{n}{v},\C{n-1}{v}:v\in V,\;n\ge1\}
\eqd
\{\W{v}{n},H_{(v,n-1)}:v\in V,\;n\ge1\},
$$
and the proof is complete.
\end{proof}

\begin{theorem}[Convergence of interacting urn schemes with drawing-based interaction]
\label{thm:convergence_drawing_interaction_urn}
Let
$
\{\C{n}{v}:v\in V,\;n\ge0\}
$
be the interacting urn scheme defined in
\Cref{subsec:drawing_interaction}. Assume that
$$
\sup_{v\in V}T_0^{(v)}<\infty,
$$
and that the reinforcement matrix $R$ is balanced, irreducible, and aperiodic. Then, for every $v\in V$,
$$
\frac{\C{n}{v}}{T_n^{(v)}}
\xrightarrow{a.s.}
\pi,
\hspace{0.9cm}
\text{as }n\to\infty,
$$
where $\pi$ is the unique probability vector satisfying
$$
\pi R=r\pi.
$$
\end{theorem}

\begin{proof}
Let
$
S:=R/r.
$
First observe that
$$
\frac{\C{n}{v}-\C{0}{v}}{T_n^{(v)}}
=
\frac{nr}{T_n^{(v)}}
\sum_{i=1}^k
S_i
\left(
\frac1n
\sum_{p=1}^n
\bb1_{\{\Z{p}{v}=i\}}
\right).
$$
Hence, by the representation \Cref{thm:reprentation_3}, it is enough to prove that
$$
\frac1n
\sum_{p=1}^n
\bb1_{\{\W{v}{p}=i\}}
\xrightarrow{a.s.}
\pi_i,
\hspace{0.9cm}
\forall v\in V,\;
i=1,\ldots,k,
$$
where
$
\{\W{v}{n}:v\in V,\;n\ge1\}
$
is the branching Markov process defined in
\Cref{subsubsec:asasd}.

For every $t>0$, let
$$
\bF_w^{\,t}
\sim
\operatorname{IRRF}_3^w
\left(
\G,
\{t\},
\{w_{\e}:\e\in E\}
\right).
$$
The proof now proceeds exactly along the same lines as that of
\Cref{thm:main_result}, provided \Cref{lem:admissible_comparison,lem:comparison_with_IRRF2,lem:comparison_with_hoppe_tree,lem:power_distance_comparison,lem:covariance_bound}
remain valid when the forest
$$
\bF
\sim
\operatorname{IRRF}_3^w
\left(
\G,
\{c_v:v\in V\},
\{w_{\e}:\e\in E\}
\right),
$$
with root weights
$
c_v=T_0^{(v)}/r,
$
is used in place of
$\operatorname{IRRF}_1$ or
$\operatorname{IRRF}_2$, and
$\bF^t$ is replaced by
$\bF_w^t$.

Since \Cref{lem:admissible_comparison} holds for $\text{IRRF}_3$,
it can be checked that the proof goes through when
$1/d_v^{\text{in}}$ is replaced by corresponding weight $w_{(u,v)}$.
Therefore, \Cref{lem:admissible_comparison} continues to hold with
$\bF^c$ and $\bF^1$ replaced by $\bF_w^c$ and $\bF_w^1$, respectively.
Since \Cref{lem:comparison_with_IRRF2} follows directly from
\Cref{lem:admissible_comparison}, it also remains valid in the present
setting.

Next consider \Cref{lem:comparison_with_hoppe_tree}.
The only modification required in its proof concerns the vertex components
$\{V_l^X\}$
and
$\{V_l^Y\}$
of the backward ancestral chains
$\{X_l\}$
and
$\{Y_l\}$, respectively.
These are now (weighted) in-neighboring random walks with transition probabilities
$$
\Q\!\left(
V_l^X=u
\,\middle|\,
V_{l-1}^X=g
\right)
=
\Q\!\left(
V_l^Y=u
\,\middle|\,
V_{l-1}^Y=g
\right)
=
w_{(u,g)}
\bb1_{\{u\in N_g\}}.
$$
On the other hand, the distributions of the time components
$\{T_l^X\}$
and
$\{T_l^Y\}$
remain unchanged.
Proceeding exactly as in the proof of \Cref{lem:comparison_with_hoppe_tree}
shows that the same conclusion holds for
$\bF_w^t$.
Consequently, \Cref{cor:properties_depth_LCA}
also remains valid for the forest
$\bF$.

At this point, it is clear that the above arguments do not depend on the
particular choice
$
w_{(u,v)}=1/d_v^{\text{in}},
$
and remain valid for arbitrary edge weights.
For \Cref{lem:power_distance_comparison},
the only change occurs in the connection probabilities appearing in the proof.
For instance, the expression for
$\Q(A_\alpha)$
changes accordingly, but all the estimates
\eqref{bound1},
\eqref{bound2},
and
\eqref{bound3}
continue to hold after replacing
$1/d_v^{\text{in}}$
by
$w_{(u,v)}$.
Hence \Cref{lem:power_distance_comparison}
also remains valid for
$\bF$.

Finally, notice that \Cref{lem:covariance_bound}
is not specific to
$\operatorname{IRRF}_1$
or
$\operatorname{IRRF}_2$ or $\operatorname{IRRF}_2$.
Its proof uses only the depth and distance estimates for the underlying recursive forest and the geometric ergodicity of the Markov chain with transition matrix
$S$. Hence \Cref{lem:covariance_bound}
holds for
$\bF$
as well.
The proof is now completed exactly as the proof of
\Cref{thm:main_result}. 
\end{proof}

\begin{remark}
\textnormal{\Cref{thm:convergence_drawing_interaction_urn} extends Theorem~4.1 of
\cite{aletti2017interacting} from finite complete graphs to arbitrary infinite
locally finite directed graphs. One can also formulate a measure-valued analog
of this interacting urn scheme, analogous to the model introduced in \Cref{subsec:IMVU_definition}, and establish a corresponding convergence
result (\Cref{thm:IMVU_convergence}) by using the tightness property of the depth of lowest common ancestor in $\bF$ (\Cref{cor:properties_depth_LCA}).}
\end{remark}

\subsection{Elephant random walks with graph-based coupled memory}
\label{subsec:ERW_coupled_memory}

\newcommand{\X}[2]{X_{#1}^{(#2)}}
\newcommand{\yi}[2]{\xi_{#1}^{(#2)}}
\renewcommand{\D}[2]{d_{#1}^{(#2)}}
\newcommand{\s}[2]{S_{#1}^{(#2)}}
\renewcommand{\us}[2]{u_{#1}^{(#2)}}
\renewcommand{\st}[2]{\hat{S}_{#1}^{(#2)}}
\newcommand{\Xt}[2]{\hat{X}_{#1}^{(#2)}}
\newcommand{\dl}[2]{\delta_{#1}^{(#2)}}
The elephant random walk is a simple example of a random walk with memory. It
was first introduced in \cite{Trimper} and was later shown to be closely
connected to randomized urn models (see \cite{Bertoin}). Variants involving
multiple elephants with shared memory have also been studied in the literature.
In this section, we show how such models are connected to
interacting urn schemes.

Suppose that $\G=(V,\E)$ is a locally finite graph as in the previous sections,
with
$
d_v^{\text{in}}>0$
for all $v\in V
$
and let
$
\{w_{\e}:\e\in E\}
$
be the collection of edge weights introduced in
\Cref{subsec:drawing_interaction}.

Let
$
\{p_v:v\in V\}$ and
$\{q_v:v\in V\}
$
be collections of real numbers in $[0,1]$. Consider a collection of independent
random variables
$$
\left\{
\sigma^{(v)},
\yi{n}{v},
\D{n}{v},
\us{n}{v}
:
n\ge1,\;v\in V
\right\},
$$
such that
\begin{itemize}
    \item[(i)] $\sigma^{(v)}\sim\operatorname{Rad}(q_v)$ for every $v\in V$;
    \item[(ii)] $\yi{n}{v}\sim\operatorname{Rad}(p_v)$ for every
    $v\in V$ and $n\ge1$;
    \item[(iii)] $\D{n}{v}\sim\operatorname{Unif}(\{1,\ldots,n\})$
    for every $v\in V$ and $n\ge1$;
    \item[(iv)]
    $
    \Q\!\left(\us{n}{v}=g\right)
    =
    w_{(g,v)}\bb1_{\{g\in N_v\}}
    $ for every $v\in V$ and $n\ge1$.
\end{itemize}
Here $\operatorname{Rad}(p)$ denotes the Rademacher distribution taking
the values $+1$ and $-1$ with probabilities $p$ and $1-p$, respectively,
while $\operatorname{Unif}(\{1,\ldots,n\})$ denotes the discrete uniform
distribution on $\{1,\ldots,n\}$.

To each vertex $v\in V$, associate a random walk
$
\{\s{n}{v}:n\ge0\},
$
whose increments
$
\{\X{n}{v}:n\ge1\}
$
take values in $\{-1,+1\}$ and are defined recursively as follows.
Set
$$
\X{1}{v}=\sigma^{(v)},
\hspace{0.9cm}
\forall v\in V.
$$
For $n\ge1$, after defining
$
\{\X{j}{v}:v\in V,\;j\le n\},
$
set
$$
\X{n+1}{v}
=
\yi{n}{v}\,
\X{\D{n}{v}}{\us{n}{v}},
\hspace{0.9cm}
\forall v\in V.
$$
Consequently,
$$
\s{0}{v}=0,
\hspace{0.9cm}
\s{n}{v}
=
\sum_{j=1}^n\X{j}{v},
\hspace{0.9cm}
\forall v\in V,\;n\ge1.
$$

The collection of random walks
$
\{\s{n}{v}:v\in V,\;n\ge0\}
$
is called the \emph{system of elephant random walks with coupled memory},
or, briefly, a \emph{system of coupled elephants}, associated with the graph
$\G$, memory parameters $\{p_v:v\in V\}$, initial step probabilities
$\{q_v:v\in V\}$, and edge weights
$\{w_{\e}:\e\in E\}$. For each $v\in V$,
$
\{\s{n}{v}:n\ge0\}
$
is called the \emph{elephant random walk} associated with the vertex $v$.

\begin{remark}
\textnormal{For a finite complete graph $\G$, this model was first introduced in
\cite{Marquioni2019ephantrandomwalkscoupled}. For arbitrary finite graphs
with uniform edge weights
$
w_{(u,v)}=1/d_v^{\text{in}},
$
the first- and second-order asymptotic properties were studied in
\cite{das2024elephantrandomwalksgraph}. More recently, the two-vertex
complete graph with general edge weights was studied in
\cite{Aguetch_Qin}.}
\end{remark}

\subsubsection{Interacting urn representation of coupled elephants}

Assume now that
$
p_v=p$ for all $v\in V$.
Let
$$
R:=
\begin{pmatrix}
p & 1-p\\
1-p & p
\end{pmatrix}.
$$
We associate with the system of coupled elephants an interacting urn
scheme as follows. Let
$$
\C{0}{v}
:=
\frac12
\left(
1+(2p-1)\X{1}{v},
1-(2p-1)\X{1}{v}
\right),
\hspace{0.9cm}
v\in V,
$$
and, for $n\ge1$, define
$$
\Z{n}{v}
:=
\frac{3-\X{n+1}{v}}{2}.
$$
Finally, set
$$
\C{n}{v}
=
\C{0}{v}
+
\sum_{j=1}^n
R_{\Z{j}{v}},
\hspace{0.9cm}
\forall n\ge1.
$$

For $i=1,2$, we then have
$$
\C{n}{v}(i)
=
\sum_{j=1}^{n+1}
\frac{
1+(-1)^{i+1}(2p-1)\X{j}{v}
}{2}.
$$
In particular,
$$
T_n^{(v)}
=
\C{n}{v}(1)+\C{n}{v}(2)
=
n+1.
$$

It is straightforward to verify that, conditioned on
$
\{\C{0}{v}:v\in V\},
$
or, equivalently, under the regular conditional probability given
$
\sigma(\{\C{0}{v}:v\in V\}),
$
the process
$
\{\C{n}{v}:v\in V,\;n\ge0\}
$
is an interacting urn scheme with drawing-based interaction, with initial
configurations
$
\{\C{0}{v}:v\in V\},
$
reinforcement matrix $R$, and edge weights
$
\{w_{\e}:\e\in E\},
$
as introduced in
\Cref{subsec:drawing_interaction}. Here $k=2$, and
$
\{\Z{n}{v}:v\in V,\;n\ge1\}
$
is the corresponding color sequence.

Notice that, for every $p\in(0,1)$, the matrix $R$ is irreducible and aperiodic.
Therefore, by
\Cref{thm:convergence_drawing_interaction_urn},
$$
\frac{\C{n}{v}}{n+1}
\xrightarrow{a.s.}
\left(\frac12,\frac12\right),
\hspace{0.9cm}
v\in V,
$$
under the regular conditional probability given
$
\sigma(\{\C{0}{v}:v\in V\}).
$
Since the limiting vector does not depend on the initial configuration,
we obtain the following result.

\begin{theorem}[Strong convergence of coupled elephants]
\label{thm:ERWG_convergence}
Let $\G=(V,\E)$ be a locally finite graph satisfying
$
d_v^{\text{in}}>0
$ for all $v\in V$. Suppose that
$
\{\s{n}{v}:v\in V,\;n\ge0\}
$
is a system of coupled elephants associated with $\G$, with a common
memory parameter $p\in(0,1)$, initial step probabilities
$\{q_v:v\in V\}$, and edge weights
$\{w_{\e}:\e\in E\}$. Then, for every $v\in V$,
$$
\frac{\s{n}{v}}{n}
\xrightarrow{a.s.}
0,
\hspace{0.9cm}
\text{as}~n\to\infty.
$$
\end{theorem}

\begin{remark}
\textnormal{\Cref{thm:ERWG_convergence} extends the almost sure convergence result
of \cite{das2024elephantrandomwalksgraph} from finite graphs to infinite
locally finite graphs, in the case where all elephants share a common
memory parameter $p\in(0,1)$.}
\end{remark}

\section{Proof of the results related to Hoppe trees}\label{sec:proofs_Hoppe_trees}
\noindent\textbf{Proof of \Cref{lem:Hoppe_depth}.} 
It is shown in the proof of Theorem 2.1 of \cite{Hoppetree} that
$$
|n|_{\mathcal T^t}
\stackrel{d}{=}
1+\sum_{j=1}^{n-1}X_{j,n},
$$
where $X_{1,n},\ldots,X_{n-1,n}$ are independent Bernoulli random variables satisfying
$$
\Q(X_{j,n}=1)=\frac1{t+j}.
$$
Therefore,
$$
\mathbb E\left[|n|_{\mathcal T^t}\right]
=
1+\sum_{j=1}^{n-1}\frac1{t+j}
\sim \log n,
$$
and
$$
\text{Var}\left(|n|_{\mathcal T^t}\right)
=
\sum_{j=1}^{n-1}\text{Var}(X_{j,n})
\le
\sum_{j=1}^{n-1}\mathbb E[X_{j,n}]
=
\mathbb E\left[|n|_{\mathcal T^t}\right].
$$
Hence, by Chebyshev's inequality,
$$
\frac{
|n|_{\mathcal T^t}
-
\mathbb E\left[|n|_{\mathcal T^t}\right]
}{
\mathbb E\left[|n|_{\mathcal T^t}\right]
}
\xrightarrow{\mathbb P}0.
$$
Since
$$
\mathbb E\left[|n|_{\mathcal T^t}\right]
\sim
\log n,
$$
the result follows. \qed\\

\noindent\textbf{Proof of \Cref{lem:tightness_LCA_Hoppe_tree}.} Suppose that the graph $\G$ consists of a single vertex $v$ with a self-loop,
that is,
$
V=\{v\}$ and
$\E=\{(v,v)\}.
$
Let
$$
\bF^t\sim\operatorname{IRRF}_3(\G,t).
$$
By \Cref{cor:properties_depth_LCA}, we have
$$
\limsup_{n>m\to\infty}
\Q\!\left(
\left|(v,m)\wedge(v,n)\right|_{\bF^t}>K
\right)
\leq
\frac{1}{2^{K-1}}.
$$

Consider the backward ancestral chains $\{X_l\}$ and $\{Y_l\}$ in $\bF^t$,
started from $(v,m)$ and $(v,n)$, respectively. Since the vertex parts of both
chains remain fixed at $v$, while their time parts have the same joint
distribution as the corresponding backward ancestral chains in a Hoppe tree
$\T{}{t}$ started from $m$ and $n$, respectively, we immediately have
$$
\left\{
\left|(v,m)\wedge(v,n)\right|_{\bF^t}
:\ n>m\ge1
\right\}
\eqd
\left\{
\left|m\wedge n\right|_{\T{}{t}}
:\ n>m\ge1
\right\}
$$
which completes the proof.\qed

\bibliographystyle{alpha}
\bibliography{IUIG}
\end{document}